\documentclass{birkjour}
\usepackage{graphics,bbm,,amssymb,amsmath,amsfonts,dsfont}

\usepackage[T1]{fontenc}             

\usepackage{graphics,color} 

\usepackage{epsfig}

\usepackage{hyperref}

\usepackage{appendix}
\usepackage{xspace}

\newtheorem{theorem}{Theorem}

\newtheorem{proposition}[theorem]{Proposition}

\newtheorem{corollary}[theorem]{Corollary}

\numberwithin{theorem}{section}
\numberwithin{claim}{section}

\newcommand{\thmref}[1]{Theorem~\ref{thm:#1}} 
\newcommand{\propref}[1]{Proposition~\ref{prop:#1}} 
\newcommand{\corref}[1]{Corollary~\ref{cor:#1}} 
\newcommand{\secref}[1]{Section~\ref{sec:#1}} 
\newcommand{\figref}[1]{Figure~\ref{fig:#1}} 

\def\R{\mathbb{R}} 
\def\N{\mathbb{N}} 
\def\Xbb{\mathbb{X}}

\providecommand{\norm}[1]{\left \lVert #1 \right\rVert}

\newcommand{\E}{\mathbb{E}} 
\newcommand{\Ex}[1]{\E\croch{#1}}

\renewcommand{\Pr}[1]{\mathbb{P}\left(#1\right)} 

\newcommand{\petito}{\mathrm{o}}

\def\sA{\mathcal{A}}\def\sC{\mathcal{C}}
\def\sE{\mathcal{E}}
\def\sH{\mathcal{H}}\def\sI{\mathcal{I}}
\def\sJ{\mathcal{J}}\def\sK{\mathcal{K}}

\def\sS{\mathcal{S}}
\def\sW{\mathcal{W}}\def\sX{\mathcal{X}}

\def\vphi{\varphi}

\newcommand{\paren}[1]{\left( \left. #1 \right. \right)} 
\newcommand{\croch}[1]{\left[ \left. #1 \right. \right]} 
\newcommand{\set}[1]{\left\{ \left. #1 \right. \right\}}
\newcommand{\absj}[1]{\left\lvert #1 \right\rvert} 

\newcommand{\psh}[2]{\ensuremath{\langle #1,#2\rangle}}

\newcommand{\egaldef}{:=} 

\newcommand{\telque}{\, \mbox{ s.t. } \,} 

\newcommand{\PEInf}[1]{\left\lfloor#1\right\rfloor} 

\def\bayes{s}

\def\kK{k\in\sK}
\def\kh{{\widehat{k}}}

\DeclareMathOperator{\pen}{pen}

\newcommand{\penid}{\pen_{\mathrm{id}}} 
\newcommand{\penmin}{\pen_{\mathrm{min}}} 
\newcommand{\penopt}{\pen_{\mathrm{opt}}} 

\def\ERM{\widehat{s}}

\def\Boule{\mathbb{B}}

\begin{document}

%
%
%
%
%
%
%
%
%

\title[Optimal kernel selection]
 {Optimal kernel selection for density estimation}

\author{M. Lerasle}

\address{%
Univ. Nice Sophia Antipolis LJAD CNRS UMR 7351\\
06100 Nice France}

\email{mlerasle@unice.fr}

\thanks{This research was partly supported by the french Agence Nationale de la Recherche (ANR 2011 BS01 010 01 projet Calibration).}
\author{N. Magalh\~aes}
\address{INRIA, Select Project\\
Univ. Paris-Sud 11\\
Departement de Mathematiques d'Orsay\\
91405 Orsay Cedex - France}
\email{nelo.moltermagalhaes@gmail.com}
\author{P. Reynaud-Bouret}
\address{Univ. Nice Sophia Antipolis LJAD CNRS UMR 7351\\
06100 Nice France}
\email{Patricia.REYNAUD-BOURET@unice.fr}



\begin{abstract}
We provide new general kernel selection rules thanks to penalized least-squares criteria. We derive optimal oracle inequalities using adequate concentration tools. We also investigate the problem of minimal penalty as described in \cite{Bir_Mas:2007}. 
\end{abstract}

\maketitle
\section{Introduction}

Concentration inequalities are central in the analysis of adaptive nonparametric statistics.  They lead to sharp penalized criteria for model selection \cite{Massart2007}, to select bandwidths and even approximation kernels for Parzen's estimators in high dimension \cite{GL2011}, to aggregate estimators \cite{RigTsy2007} and to properly calibrate thresholds \cite{DJKP1996}. 

\noindent In the present work, we are interested in the selection of a general kernel estimator based on a least-squares density estimation approach. 
The problem has been considered in $L^1$-loss by Devroye and Lugosi \cite{Devroye-Lugosi2001}. Other methods combining log-likelihood and roughness/smoothness penalties have also been proposed in  \cite{Egg_LaR:1999,Egg_LaR:1999b,Egg_LaR:2001}. However these estimators  are usually quite difficult to compute in practice. We propose here to minimize penalized least-squares criteria and obtain from them more easily computable estimators.
Sharp concentration inequalities for U-statistics \cite{GLZ2000,Adamczak2006,HRB2003} control the variance term of the kernel estimators, whose asymptotic behavior has been precisely described, for instance in \cite{MasSwan2011,MasSwan2015,DehOuadah2013}.
We derive from these bounds (see \propref{ConcRisk}) a penalization method to select a kernel which satisfies an asymptotically optimal oracle inequality, i.e. with leading constant asymptotically equal to $1$. 

\noindent In the spirit of \cite{GinNickl2009}, we use an extended definition of kernels that allows to deal simultaneously with classical collections of estimators as projection estimators, weighted projection estimators, or Parzen's estimators.
This method can be used for example to select an optimal model in model selection (in accordance with \cite{Massart2007}) or to select an optimal  bandwidth together with an optimal approximation kernel among a finite collection of Parzen's estimators. In this sense, our method deals, in particular, with the same problem as that of Goldenshluger and Lepski \cite{GL2011} and we establish in this framework that a leading constant 1 in the oracle inequality is indeed possible.

\noindent Another main consequence of concentration inequalities is to prove the existence of  a minimal level of penalty, under which no oracle inequalities can hold. Birg\'e and Massart shed light on this phenomenon in a Gaussian setting for model selection  \cite{Bir_Mas:2007}. Moreover in this setting, they  prove that the optimal penalty is twice the minimal one. In addition, there is a sharp phase transition in the dimension of the selected models leading to an estimate of the optimal penalty in their case (which is known up to a multiplicative constant). Indeed, starting from the idea that in many models the optimal penalty is twice the minimal one (this is the slope heuristic), Arlot and Massart \cite{Arl_Mas:2009} propose to detect the minimal penalty by the phase transition and to apply the "$\times 2$" rule (this is the slope algorithm). They prove that this algorithm works at least in some regression settings.

\noindent In the present work, we also show that  minimal penalties exist in the density estimation setting. In particular, we exhibit a sharp "phase transition" of the behavior of the selected estimator around this minimal penalty. 
The analysis of this last result is not standard however. First, the "slope heuristic" of \cite{Bir_Mas:2007} only holds in particular cases as the selection of projection estimators, see also  \cite{Le2012}. As in the selection of a linear estimator in a regression setting \cite{Arl_Bac:2009}, the heuristic can sometimes be corrected: for example for the selection of a bandwidth when the approximation kernel is fixed. In general since there is no simple relation between the minimal penalty and the optimal one, the slope algorithm of \cite{Arl_Mas:2009} shall only be used with care for kernel selection.
Surprisingly our work reveals that the minimal penalty can be negative. In this case, minimizing an unpenalized criterion leads to oracle estimators. To our knowledge, such phenomenon has only been noticed previously in a very particular classification setting \cite{FrTul2006}. We illustrate all of these different behaviors by means of a simulation study. 
 
\noindent In \secref{not}, after fixing the main notation, providing some examples and defining the framework, we explain our goal, describe what we mean by an \textit{oracle inequality} and state the exponential inequalities that we shall need. Then we derive optimal penalties in \secref{opt} and study the problem of minimal penalties in \secref{MinPen}. All of these results are illustrated for our three main examples : projection kernels, approximation kernels and weighted projection kernels. In \secref{simu}, some simulations are performed in the approximation kernel case. The main proofs are detailed in \secref{proofs} and technical results are discussed in the appendix.

%

\section{Kernel selection for least-squares density estimation \label{sec:not}}

\subsection{Setting}
Let $X,Y,X_1,\ldots,X_n$ denote i.i.d. random variables taking values in the measurable space $(\Xbb,\sX,\mu)$, with common distribution $P$. Assume $P$ has density $\bayes$ with respect to $\mu$ and $\bayes$ is uniformly bounded. Hence, $\bayes$ belongs to $L^2$, where, for any $p\ge 1$,
\vspace{-0.2cm}
\[L^p\egaldef\set{t:\Xbb\to\R,\telque  \norm{t}_p^p \egaldef \int\absj{t}^pd\mu<\infty}\enspace.\]
\vspace{-0.4cm}

\noindent Moreover, $\norm{\cdot}=\norm{\cdot}_2$ and $\psh{\cdot}{\cdot}$ denote respectively the $L^2$-norm and the associated inner product and $\norm{\cdot}_\infty$ is the supremum norm. We  systematically use  $x\vee y$ and $x\wedge y$ for $\max(x,y)$ and $\min(x,y)$ respectively, and denote  $|A|$ the cardinality of the set $A$. 
Recall that $x_+=x\vee 0$ and, for any $y\in \R^+$, $\PEInf{y}=\sup\{n\in \N \telque n\le y\}$. 

\smallskip

\noindent Let $\set{k}_{\kK}$ denote a collection of symmetric functions $k:\Xbb^2\to\R$ indexed by some given finite set $\mathcal{K}$ such that 
\[\sup_{x\in\Xbb}~\int_\Xbb k(x,y)^2d\mu(y)~\vee \sup_{(x,y)\in\Xbb^2}\absj{k(x,y)}<\infty\enspace.\]
A function $k$ satisfying these assumptions is called a \textit{kernel}, in the sequel. 
A kernel $k$ is associated with an estimator $\ERM_k$ of $\bayes$ defined for any $x\in\Xbb$ by
\vspace{-0.1cm}
\begin{align*}
 \ERM_k(x)&\egaldef \frac1n\sum_{i=1}^nk(X_i,x) \enspace.
\end{align*}
\vspace{-0.3cm}

\noindent Our aim is to select a ``good'' $\ERM_{\hat{k}}$ in the family $\{\ERM_k,k \in \mathcal{K}\}$.
 Our results are expressed in terms of a constant $\Gamma\ge 1$ such that for all $\kK$,
\begin{equation}\label{eq:CdtForGamma}
\sup_{x\in\Xbb}~\int_\Xbb k(x,y)^2d\mu(y)~\vee \sup_{(x,y)\in\Xbb^2}\absj{k(x,y)}\le\Gamma n\enspace.
\end{equation}
\vspace{-0.3cm}

\noindent 
This condition plays the same role as 
\(\int |k(x,y)| s(y) d\mu(y) < \infty\), the milder condition 
used in \cite{Devroye-Lugosi2001} when working with $L^1$-losses.
Before describing the method, let us give three examples of such estimators that are used for density estimation, and see how they can naturally be associated to some kernels. \secref{ProofsKernels} of the appendix gives the computations leading to the corresponding $\Gamma$'s.

\smallskip

\paragraph{Example 1: Projection estimators.}
Projection estimators are among the most classical density estimators.
Given a linear subspace $S\subset L^2$, the projection estimator on $S$ is defined by 
\[\ERM_S=\arg\min_{t\in S}\set{\norm{t}^2-\frac2n\sum_{i=1}^nt(X_i)}\enspace.\]
Let $\mathcal{S}$ be a family of linear subspaces $S$ of $L^2$. For any $S\in\sS$, let $(\vphi_\ell)_{\ell\in\sI_S}$ denote an orthonormal basis of $S$. The projection estimator $\ERM_S$ can be computed and is equal to 
\[\ERM_S=\sum_{\ell\in\sI_S}\paren{\frac1n\sum_{i=1}^n\vphi_\ell(X_i)}\vphi_\ell\enspace.\]
It is therefore easy to see that it is the estimator associated to the \textit{projection kernel} $k_S$ defined for any $x$ and $y$ in $\Xbb$ by
\[k_S(x,y)\egaldef \sum_{\ell\in\sI_S}\vphi_\ell(x)\vphi_\ell(y)\enspace.\]
\vspace{-0.2cm}

\noindent Notice that $k_S$ actually depends on the basis $(\vphi_\ell)_{\ell\in\sI_S}$ even if $\ERM_S$ does not. In the sequel, we always assume that some orthonormal basis $(\vphi_\ell)_{\ell\in\sI_S}$ is given with $S$.
Given a finite collection $\sS$ of linear subspaces of $L^2$, one can choose the following constant $\Gamma$ in \eqref{eq:CdtForGamma} for the collection $(k_S)_{S\in\sS}$  
\begin{equation}\label{gammaproj}
\Gamma=1\vee \frac 1n\sup_{S\in \sS}\sup_{f\in S, \norm{f}=1} \norm{f}_\infty^2\enspace.
\end{equation}
\vspace{-0.2cm}
\paragraph{Example 2: Parzen's estimators.}
Given a bounded symmetric integrable function $K:\R\to\R$ such that $\int_\R K(u)du=1$, $K(0)>0$ and a bandwidth $h>0$, the Parzen estimator is defined  by
\vspace{-0.2cm}
\[\forall x \in \R, \quad \ERM_{K,h}(x)=\frac1{nh}\sum_{i=1}^nK\paren{\frac{x-X_i}h}\enspace.\] 
\vspace{-0.2cm}

\noindent It can also naturally be seen as a kernel estimator, associated to the function $k_{K,h}$ defined for any $x$ and $y$ in $\R$ by
\vspace{-0.1cm}
\[k_{K,h}(x,y)\egaldef \frac1hK\paren{\frac{x-y}h}\enspace.\]
\vspace{-0.2cm}

\noindent We shall call the function $k_{K,h}$ an approximation or Parzen kernel.\\
Given a finite collection of pairs $(K,h)\in\sH$, one can choose $\Gamma=1$ in \eqref{eq:CdtForGamma} if, 
\begin{equation}\label{eq:cond_gamma_Kh}
h\ge \frac{\norm{K}_\infty\norm{K}_1}{n} \quad\mbox{ for any }(K,h)\in\sH\enspace.
\end{equation}
\vspace{-0.2cm}
\paragraph{Example 3: Weighted projection estimators.}
Let $(\vphi_i)_{i=1,\ldots,p}$ denote an orthonormal system in $L^2$ and let $w=(w_i)_{i=1,\ldots,p}$ denote real numbers in $[0,1]$. The associated weighted kernel projection estimator of $\bayes$ is defined by
\[\ERM_w=\sum_{i=1}^pw_i\paren{\frac1n \sum_{j=1}^n \vphi_i(X_j)}\vphi_i\enspace.\]
These estimators are used to derive very sharp adaptive results. In particular, Pinsker's estimators are weighted kernel projection estimators (see for example \cite{Rig2006}). When $w\in\set{0,1}^p$, we recover a classical projection estimator. A weighted projection estimator is associated to the \textit{weighted projection kernel} defined for any $x$ and $y$ in $\Xbb$ by
\[k_w(x,y)\egaldef \sum_{i=1}^p w_i\vphi_i(x)\vphi_i(y) \enspace.\]
Given any finite collection $\mathcal{W}$ of weights, one can choose in \eqref{eq:CdtForGamma}
\begin{equation}\label{gammaweight}
\Gamma=1\vee \paren{\frac1n\sup_{x\in \Xbb}\sum_{i=1}^p\vphi_i(x)^2}\enspace.
\end{equation}

\subsection{Oracle inequalities and penalized criterion}
The goal is to estimate $\bayes$ in the best possible way using a finite collection of kernel estimators $(\ERM_k)_{\kK}$. In other words, the purpose is to select among  $(\ERM_k)_{\kK}$ an estimator $\ERM_\kh$ from the data such that  $\norm{\ERM_\kh-\bayes}^2$  is as close as possible to $\inf_{\kK}\norm{\ERM_k-\bayes}^2$. More precisely our aim is to select $\kh$ such that, with high probability,
\begin{equation}\label{eq:oracledef}
\norm{\ERM_{\kh}-\bayes}^2\le C_n\inf_{\kK}\norm{\ERM_k-\bayes}^2+R_n\enspace,
\end{equation}
where $C_n\geq 1$ is the leading constant and $R_n>0$ is usually a remaining term. In this case, $\ERM_{\kh}$  is said to satisfy an \textit{oracle inequality}, as long as $R_n$ is small compared to 
$\inf_{\kK}\norm{\ERM_k-\bayes}^2$ and $C_n$ is a bounded sequence. This means that the selected estimator does as well as the best estimator in the family up to some multiplicative constant. The best case one can expect is to get $C_n$ close to 1. This is why, when $C_n\to_{n\to\infty}1$, the corresponding oracle  inequality is called \textit{asymptotically optimal}. To do so, we study minimizers of \textit{penalized least-squares criteria}. Note that in our three examples choosing $\kh \in \mathcal{K}$ amounts to choosing the smoothing parameter, that is respectively to choosing $\widehat{S}\in\sS$, $(\widehat{K},\widehat{h})\in\sH$ or $\widehat{w}\in\sW$.


\noindent Let $P_n$ denote the empirical measure, that is, for any real valued function $t$, 
\vspace{-0.1cm}
\[P_n(t)\egaldef \frac1n\sum_{i=1}^nt(X_i)\enspace.\]
\vspace{-0.3cm}

\noindent For any $t\in L^2$, let also $P(t)\egaldef \int_{\Xbb}t(x)\bayes(x)d\mu(x)\enspace.$\\
The \textit{least-squares contrast} is defined, for any $t\in L^2$, by
$$\gamma(t)\egaldef \norm{t}^2-2t\enspace.$$ 
Then for any given function $\pen: \sK \to \R$, the \textit{least-squares penalized criterion} is defined by
\begin{equation}\label{eq:CritPen}
\sC_{\pen}(k)\egaldef P_n\gamma(\ERM_k)+\pen(k)\enspace.
\end{equation}
Finally the selected $\kh\in \sK$ is given by any minimizer of $\sC_{\pen}(k)$, that is, 
\begin{equation}\label{eq:defkhat}
\kh\in\arg\min_{\kK}\set{\sC_{\pen}(k)}\enspace.
\end{equation}
\vspace{-0.3cm}

\noindent
As $P\gamma(t)=\norm{t-\bayes}^2-\norm{\bayes}^2$,  it is equivalent to minimize $\norm{\ERM_k-\bayes}^2$ or $P\gamma(\ERM_k)$. As our goal is to select $\ERM_\kh$ satisfying an oracle inequality, an ideal penalty $\penid$ should satisfy $\sC_{\penid}(k)=P\gamma(\ERM_k)$, i.e. criterion \eqref{eq:CritPen} with
\[\penid(k)\egaldef (P-P_n)\gamma(\ERM_k)=2(P_n-P)(\ERM_k)\enspace.\] 
To identify the main quantities of interest, let us introduce some notation and  develop $\penid(k)$. 
For all $\kK$, let
\begin{align*}
 \bayes_k(x)&\egaldef \int_\Xbb k(y,x)\bayes(y)d\mu(y)=\Ex{k(X,x)},\qquad \forall x\in\Xbb\enspace,
\end{align*}
 and 
\[U_k\egaldef \sum_{i\neq j=1}^n\paren{k(X_i,X_j)-\bayes_k(X_i)-\bayes_k(X_j)+\Ex{k(X,Y)}}\enspace.\]
Because those quantities are fundamental in the sequel, let us also define $\Theta_k(x)=A_k(x,x)$ where for $(x,y)\in \Xbb^2$
\vspace{-0.1cm}
\begin{equation}\label{eq:DefAk}
A_k(x,y)\egaldef     \int_{\Xbb}k(x,z)k(z,y)d\mu(z)\enspace.
\end{equation}
\vspace{-0.3cm}

\noindent Denoting $$\mbox{for all } x\in\Xbb,\quad \chi_k(x)=k(x,x)\enspace,$$
the ideal penalty is then equal to
\vspace{-0.1cm}
\begin{align}
\notag &\penid(k)=2(P_n-P)(\ERM_k-\bayes_k)+2(P_n-P)\bayes_k\\
\label{eq:penid} &=2\paren{\frac{P\chi_k-P\bayes_k}{n}+\frac{(P_n-P)\chi_k}n+\frac{U_k}{n^2}+\paren{1-\frac2n}(P_n-P)\bayes_k}\enspace.
\end{align}
\vspace{-0.3cm}

\noindent 
The main point is that by using concentration inequalities, 
we obtain:
\vspace{-0.1cm}
$$\penid(k) \simeq 2\paren{\frac{P\chi_k-P\bayes_k}{n}}\enspace.$$
\vspace{-0.3cm}

\noindent 
The term $P\bayes_k/n$ depends on $\bayes$ which is unknown. Fortunately, it can be easily controlled as detailed in the sequel. Therefore one can hope that the choice
\vspace{-0.1cm}
$$\pen(k)=2\frac{P\chi_k}n$$
\vspace{-0.3cm}

\noindent
is convenient. In general, this choice  still depends on the unknown density $\bayes$ but it can be easily estimated in a data-driven way by
\vspace{-0.1cm}
$$\pen(k)=2\frac{P_n\chi_k}n\enspace.$$
\vspace{-0.3cm}

\noindent
The goal of \secref{opt} is to prove this heuristic and to show that $2P\chi_k/n$ and $2P_n\chi_k/n$ are optimal choices for the penalty, that is, they lead to an asymptotically optimal oracle inequality.

\subsection{Concentration tools}\label{sec:tools}
To derive sharp oracle inequalities, we only need two fundamental concentration tools, namely a weak Bernstein's inequality and the concentration bounds for degenerate U-statistics of order two. We cite them here under their most suitable form for our purpose.

\paragraph{A weak Bernstein's inequality.}~

\vspace{-0.2cm}
\begin{proposition}\label{prop:Bern}
For any bounded real valued function $f$ and any $X_1,\ldots,X_n$ i.i.d. with distribution $P$, for any $u>0$,
$$\Pr{(P_n-P) f \geq \sqrt{\frac{2 P\paren{f^2} u}n}+\frac{\norm{f}_\infty u}{3n}}\leq \exp(-u)\enspace.$$
\end{proposition}
The proof is straightforward and can be derived from either Bennett's or Bernstein's inequality \cite{BLM2013}.

\paragraph{Concentration of degenerate U-statistics of order 2.}~

\vspace{-0.2cm}
\begin{proposition}\label{prop:ConcUstat}
Let $X,X_1,\ldots X_n$ be i.i.d. random variables defined on a Polish space $\Xbb$ equipped with its Borel $\sigma$-algebra and let $\paren{f_{i,j}}_{1\leq i\not=j\le n}$ denote bounded real valued symmetric measurable functions defined on $\Xbb^2$, such that for any $i\not = j$, $f_{i,j}=f_{j,i}$ and
\begin{equation}\label{eq:NonDegenerate}
\forall ~i,j \mbox{ s.t. }  1\leq i\neq j\le n,\qquad \Ex{f_{i,j}(x,X)}=0\qquad \mbox{for\;a.e.\;}x \mbox{ in }\Xbb\enspace.
\end{equation}
Let $U$ be the following totally degenerate $U$-statistic of order $2$,
\vspace{-0.1cm}
\[U=\sum_{1\leq i\neq j\le n}f_{i,j}(X_i,X_j)\enspace.\]
\vspace{-0.3cm}

\noindent Let $A$ be an upper bound of $\absj{f_{i,j}(x,y)}$ for any $i,j,x,y$ and 
\vspace{-0.2cm}
\begin{align*}
 B^2&=\max\paren{\sup_{i, x\in\Xbb}\sum_{j=1}^{i-1}\Ex{f_{i,j}(x,X_j)^2},\sup_{j,t\in\Xbb}\sum_{i=j+1}^{n}\Ex{f_{i,j}(X_i,t)^2}}\\
 C^2&=\sum_{1\le i\neq j\le n}\Ex{f_{i,j}(X_i,X_j)^2}\\
 D&=\sup_{(a,b)\in\sA}\Ex{\sum_{1\le i<j\le n}f_{i,j}(X_i,X_j)a_i(X_i)b_j(X_j)}\enspace,
\end{align*}
\vspace{-0.3cm}

\noindent
where $\displaystyle\sA=\set{(a,b),\telque \Ex{\sum_{i=1}^{n-1}a_i(X_i)^2}\le 1,\;\Ex{\sum_{j=2}^{n}b_j(X_j)^2}\le 1}.$\\
Then there exists some absolute constant $\kappa>0$ such that for any $u>0$, with probability larger than $1-2.7e^{-u}$,
\vspace{-0.2cm}
\begin{align*}
U\le \kappa\left(C\sqrt{u}+Du+Bu^{3/2}+Au^2\right) \enspace. 
\end{align*}
\vspace{-0.3cm}
\end{proposition}
\noindent The present result is a simplification of Theorem 3.4.8 in \cite{Gin_Nic:2015}, which provides explicit constants for any variables defined on a Polish space. It is mainly inspired by \cite{HRB2003}, where the result therein has been stated only for real variables.
This inequality actually dates back to Gin\'e,  Latala and Zinn \cite{GLZ2000}.  This result has been further generalized by  Adamczak to U-statistics of any order \cite{Adamczak2006}, though the constants are not explicit. 


\section{Optimal penalties for kernel selection\label{sec:opt}}
The main aim of this section is to show that $2P\chi_k/n$ is a theoretical optimal penalty for kernel selection, which means that if  $\pen(k)$ is close to $2P\chi_k/n$, the selected kernel $\kh$ satisfies an asymptotically optimal oracle inequality.
\vspace{-0.2cm}
\subsection{Main assumptions}\label{sec:MainAss}

To express our results in a simple form, a positive constant $\Upsilon$ is assumed to control, for any $k$ and $k'$ in $\sK$, all the following quantities.
\begin{eqnarray}
\paren{\Gamma(1+\norm{\bayes}_\infty)}\vee \sup_{\kK}\norm{\bayes_k}^2&\le &\Upsilon\enspace,\label{eq:AssOr1}\\
 P\paren{\chi_k^2} & \le &  \Upsilon n P\Theta_k\enspace,\label{eq:AssOr2}\\
 \norm{\bayes_k-\bayes_{k'}}_\infty &\le& \Upsilon\vee\sqrt{\Upsilon n}\norm{\bayes_k-\bayes_{k'}}\enspace,\label{eq:AssOr3}\\
 \Ex{A_k(X,Y)^2}&\le& \Upsilon P\Theta_k\enspace,\label{eq:AssOr4}\\
 \sup_{x\in\Xbb}~\Ex{A_k(X,x)^2}&\le& \Upsilon n \enspace,\label{eq:AssOr5}\\
 v_k^2\egaldef \sup_{t\in\Boule_k}P t^2&\le& \Upsilon\vee \sqrt{\Upsilon P\Theta_k}\enspace, \label{eq:AssVar}
\end{eqnarray}
where $\Boule_k$ is the set of functions $t$ that can be written $t(x)=\int a(z)k(z,x)d\mu(z)$ for some $a\in L^2$ with $\norm{a}\le 1$.

\noindent These assumptions may seem very intricate. They are actually  fulfilled by our three main examples under very mild conditions (see \secref{MainExamples}).

\subsection{The optimal penalty theorem}
In the sequel, $\square$ denotes a positive absolute constant whose value may change from line to line and if there are indices such as $\square_\theta$, it means that this is a positive function of $\theta$ and only $\theta$ whose value may change from line to line.

\begin{theorem}\label{thm:OptimalPenalty}
If Assumptions~\eqref{eq:AssOr1}, \eqref{eq:AssOr2}, \eqref{eq:AssOr3}, \eqref{eq:AssOr4} \eqref{eq:AssOr5}, \eqref{eq:AssVar} hold, then, for any $x\ge 1$, with probability larger than $1-\square|\sK|^2e^{-x}$, for any $\theta\in(0,1)$, any minimizer $\kh$ of the penalized criterion \eqref{eq:CritPen} satisfies the following inequality
 \begin{align}
\notag \forall \kK, \qquad (1-4\theta) \norm{\bayes-\ERM_{\kh}}^2&\le(1+4\theta)\norm{\bayes-\ERM_k}^2+\paren{\pen(k)-2\frac{P\chi_k}n}\\
\label{eq:WhateverPenalty}&-\paren{\pen\paren{\kh}-2\frac{P\chi_{\kh}}n}+\square\frac{\Upsilon x^2}{\theta n}\enspace.
\end{align}


\noindent Assume moreover that there exists $C>0$, $\delta'\geq \delta>0$ and $r\ge0$ such that for any $x\ge 1$, with probability larger than $1-Ce^{-x}$, for any $k\in\sK$,
\begin{equation}\label{eq:CondPen}
 (\delta-1)\frac{ P\Theta_k}n-\square r\frac{\Upsilon x^2}{n}\le \pen(k)-\frac{2P\chi_k}{n} \le (\delta'-1) \frac{P\Theta_k}n+\square r\frac{\Upsilon x^2}{n}\enspace.
\end{equation}
Then for all $\theta \in(0,1)$ and all $x\ge 1$, the following holds with probability at least $1-\square(C+|\sK|^2)e^{-x}$,
\begin{align*}
\frac{(\delta\wedge1)-5 \theta}{ (\delta'\vee 1)+(4+\delta')\theta} \norm{\bayes-\ERM_{\kh}}^2&\le\inf_{\kK}\norm{\bayes-\ERM_k}^2+\square \paren{r+\frac1{\theta^3}}\frac{\Upsilon x^2}{ n}\enspace.
\end{align*}
\end{theorem}
\noindent Let us make some remarks.
\begin{itemize}
\item First, this is an oracle inequality (see \eqref{eq:oracledef}) with leading constant $C_n$ and remaining term $R_n$ given by
$$C_n= \frac{(\delta'\vee 1)+(4+\delta')\theta}{ (\delta\wedge1)-5 \theta} \quad \mbox{and} \quad R_n= \square C_n (r+\theta^{-3})\frac{\Upsilon x^2}{n}\enspace,$$
as long as 
\begin{itemize}
\item $\theta$ is small enough for $C_n$ to be positive,
\item $x$ is large enough for the probability to be large and
\item $n$ is large enough for $R_n$ to be negligible. 
\end{itemize}
Typically, $r, \delta, \delta', \theta$ and $\Upsilon$ are bounded w.r.t. $n$ and $x$ has to be of the order of $\log(|\sK|\vee n)$ for the remainder to be negligible. In particular, $\sK$ may grow with $n$ as long as (i) $\log(|\sK|\vee n)^2$ remains negligible with respect to $n$ and (ii) $\Upsilon$ does not depend on $n$.
 \item If $\pen(k)=2P\chi_k/n$, that is if $\delta=\delta'=1$ and $r=C=0$ in \eqref{eq:CondPen}, the estimator $\ERM_{\kh}$ satisfies an asymptotically optimal oracle inequality i.e. $C_n\to_{n\to\infty} 1$ since $\theta$ can be chosen as close to 0 as desired. Take for instance, $\theta=(\log n)^{-1}$.
\item In general $P\chi_k$ depends on the unknown $\bayes$ and this last penalty cannot be used in practice. Fortunately, its empirical counterpart $\pen(k)=2P_n\chi_k/n$ satisfies \eqref{eq:CondPen} with $\delta=1-\theta$, $\delta'=1+\theta$, $r=1/\theta$ and $C=2|\sK|$ for any $\theta\in(0,1)$ and in particular $\theta=(\log n)^{-1}$ (see \eqref{eq:Pn-Pchik} in \propref{ConcRemainder}). Hence, the estimator $\ERM_{\kh}$ selected with this choice of penalty also satisfies an asymptotically optimal oracle inequality, by the same argument.
\item Finally, we only get an oracle inequality when $\delta>0$, that is when $\pen(k)$ is larger than $(2P\chi_k-P\Theta_k)/n$ up to some residual term. We  discuss the necessity of this condition in \secref{MinPen}.
\end{itemize}

\subsection{Main examples}\label{sec:MainExamples}
This section shows that \thmref{OptimalPenalty} can be applied in the examples. In addition, it provides the computation of $2P\chi_k/n$ in some specific cases of special interest.

\paragraph{Example 1 (continued).}
~\\
\vspace{-0.6cm}
\begin{proposition}\label{prop:AssProjKern}
Let $\set{k_S,S\in\sS}$ be a collection of projection kernels. Assumptions~\eqref{eq:AssOr1}, \eqref{eq:AssOr2}, \eqref{eq:AssOr4}, \eqref{eq:AssOr5} and \eqref{eq:AssVar} hold for any $\Upsilon\ge \Gamma(1+ \norm{\bayes}_\infty)$, where $\Gamma$ is given by \eqref{gammaproj}. In addition, Assumption~\eqref{eq:AssOr3} is satisfied under either of the following classical assumptions (see \cite[Chapter 7]{Massart2007}):
\begin{equation}\label{eq:Nested}
\forall S,S'\in \sS,\qquad  \mbox{ either } S\subset S' \mbox{ or } S'\subset S\enspace,
\end{equation}
or
\vspace{-0.2cm}
\begin{equation}\label{eq:UnifSupBound}
 \forall S\in\sS,\qquad \norm{\bayes_{k_S}}_\infty\leq \frac{\Upsilon}2\enspace.
\end{equation}
\end{proposition}

\noindent
These particular projection kernels satisfy for all $(x,y)\in \Xbb^2$
\begin{multline*}
A_{k_S}(x,y)=\int_\Xbb k_S(x,z)k_S(y,z)d\mu(z)\\
=\sum_{(i,j)\in\sI^2_S}\vphi_i(x)\vphi_j(y)\int_{\Xbb}\vphi_i(z)\vphi_j(z)d\mu(z)=k_S(x,y)\enspace.
\end{multline*}
In particular, $\Theta_{k_S}=\chi_{k_S}=\sum_{i\in\sI_S}\vphi_i^2$ and $2P\chi_{k_S}-P\Theta_{k_S}=P\chi_{k_S}$. 

\noindent Moreover, it appears that the function $\Theta_{k_S}$ is constant in some linear spaces $S$ of interest (see \cite{Le2012} for more details). Let us mention one particular case studied further on in the sequel. Suppose $\sS$ is a collection of regular histogram spaces $S$ on $\Xbb$, that is, any $S\in \sS$ is a space of piecewise constant functions on a partition $\sI_S$ of $\Xbb$ such that $\mu(i)=1/D_S$ for any $i$ in $\sI_S$. Assumption~\eqref{eq:UnifSupBound} is satisfied for this collection as soon as $\Upsilon\ge 2\norm{\bayes}_\infty$. The family $(\vphi_i)_{i\in\sI_S} $, where $\vphi_i=\sqrt{D_S}{\bf 1}_{i}$ is an orthonormal basis of $S$ and 
\vspace{-0.1cm}
\[\chi_{k_S}=\sum_{i\in\sI_S}\vphi_i^2=D_S\enspace.\]
\vspace{-0.2cm}

\noindent Hence, $P\chi_{k_S}=D_S$ and $2D_S/n$ can actually be used as a penalty to ensure that the selected estimator satisfies an asymptotically optimal oracle inequality. Moreover, in this example it is actually necessary to choose a penalty larger than $D_S/n$ to get an oracle inequality (see \cite{Le2012} or \secref{MinPen} for more details).

\paragraph{Example 2 (continued).}
~\\
\vspace{-0.6cm}
\begin{proposition}\label{prop:AssApproxKern}
Let $\set{k_{K,h},(K,h)\in\sH}$ be a collection of approximation kernels. Assumptions~\eqref{eq:AssOr1}, \eqref{eq:AssOr2}, \eqref{eq:AssOr3}, \eqref{eq:AssOr4}, \eqref{eq:AssOr5} and \eqref{eq:AssVar} hold with $\Gamma=1$, for any 
\[
\Upsilon \ge \max_{K}\set{\frac{K(0)}{\norm{K}^2}\vee \paren{1+2\norm{\bayes}_\infty\norm{K}_1^2}}\enspace,
\]
as soon as \eqref{eq:cond_gamma_Kh} is satisfied.
\end{proposition}

\noindent
These approximation kernels satisfy, for all $x\in \R$,
\begin{align*}
\chi_{k_{K,h}}(x)&=k_{K,h}(x,x)=\frac{K(0)}{h}\enspace,\\
\Theta_{k_{K,h}}(x)&=A_{k_{K,h}}(x,x)=\frac1{h^2}\int_{\R}K\paren{\frac{x-y}{h}}^2dy=\frac{\norm{K}^2}{h}\enspace. 
\end{align*}
Therefore, the optimal penalty $2P\chi_{k_{K,h}}/n=2K(0)/{(nh)}$ can be computed in practice and yields an asymptotically optimal selection criterion.
Surprisingly, the lower bound $2P\chi_{k_{K,h}}/n-P\Theta_{k_{K,h}}/n=(2K(0)-\norm{K}^2)/(nh)$ can be negative if $\norm{K}^2>2K(0)$. In this case, a minimizer of \eqref{eq:CritPen} satisfies an\ oracle inequality, even if this criterion is not penalized. This remarkable fact is illustrated in the simulation study in \secref{simu}. 

%
%

\paragraph{Example 3 (continued).}
~\\
\vspace{-0.6cm}
\begin{proposition}\label{prop:AssWeiProjKern}
Let $\set{k_w,w\in\sW}$ be a collection of weighted projection kernels. Assumption~\eqref{eq:AssOr1} is valid for $\Upsilon\ge \Gamma(1+ \norm{\bayes}_\infty)$, where $\Gamma$ is given by \eqref{gammaweight}. Moreover \eqref{eq:AssOr1} and \eqref{eq:CdtForGamma} imply \eqref{eq:AssOr2}, \eqref{eq:AssOr3}, \eqref{eq:AssOr4}, \eqref{eq:AssOr5} and \eqref{eq:AssVar}.
\end{proposition}

\noindent
For these weighted projection kernels, for all $x\in \Xbb$
\begin{align*}
\chi_{k_w}(x)&=\sum_{i=1}^pw_i\vphi_i(x)^2,\qquad\mbox{hence}\qquad P\chi_{k_w}=\sum_{i=1}^pw_iP\vphi_i^2\enspace \qquad \mbox{and}\\
\Theta_{k_w}(x)&=\sum_{i,j=1}^pw_iw_j\vphi_i\vphi_j\int_{\Xbb}\vphi_i(x)\vphi_j(x)d\mu(x)=\sum_{i=1}^pw_i^2\vphi_i(x)^2\le \chi_{k_w}(x)\enspace.
\end{align*}
In this case, the optimal penalty $2P\chi_{k_w}/n$ has to be estimated in general. However, in the following example it can still be directly computed.\\
Let $\Xbb=[0,1]$, let $\mu$ be the Lebesgue measure. Let $\vphi_0\equiv 1$ and, for any $j\ge 1$,
\[\vphi_{2j-1}(x)=\sqrt{2}\cos(2\pi j x),\qquad \vphi_{2j}(x)=\sqrt{2}\sin(2\pi j x)\enspace.\]
Consider some odd $p$ and a family of weights $\sW=\set{w_i,i=0,\ldots,p}$ such that, for any $w\in\sW$ and any $i=1,\ldots,p/2,$ $w_{2i-1}=w_{2i}=\tau_i$. In this case, the values of the functions of interest do not depend on $x$
\[\chi_{k_w}=w_0+\sum_{j=1}^{p/2}\tau_j,\qquad \Theta_{k_w}=w_0^2+\sum_{j=1}^{p/2}\tau^2_j\enspace.\]
In particular, this family includes Pinsker's and Tikhonov's weights.

%

\vspace{-0.6cm}
\section{Minimal penalties for kernel selection}\label{sec:MinPen}

The purpose of this section is to see whether the lower bound $\penmin(k)\egaldef (2P\chi_k-P\Theta_k)/n$ is sharp in \thmref{OptimalPenalty}. To do so we first need the following result which links $\norm{\bayes-\ERM_k}$ to deterministic quantities, thanks to concentration tools.

\subsection{Bias-Variance decomposition with high probability}

\begin{proposition}\label{prop:ConcRisk}
Assume $\set{k}_{\kK}$ is a finite collection of kernels satisfying Assumptions~\eqref{eq:AssOr1}, \eqref{eq:AssOr2}, \eqref{eq:AssOr3}, \eqref{eq:AssOr4} \eqref{eq:AssOr5} and \eqref{eq:AssVar}. For all $x>1$, for all $\eta$ in $(0,1]$, with probability larger than 
$1-\square |\sK|e^{-x}$
 $$ \norm{\bayes_k-\ERM_k}^2 \leq (1+\eta) \frac{P\Theta_k}{n} +\square\frac{\Upsilon x^2}{\eta n}\enspace,$$
 $$ \frac{P\Theta_k}{n} \leq (1+\eta) \norm{\bayes_k-\ERM_k}^2  +\square\frac{\Upsilon x^2}{\eta n}\enspace.$$
 
\noindent Moreover, for all $x>1$ and for all $\eta$ in $(0,1)$, with probability larger than 
 $1-\square |\sK|e^{-x}$, for all $\kK$, each of the following inequalities hold
 $$\norm{\bayes-\ERM_k}^2 \leq (1+\eta) \paren{\norm{\bayes-\bayes_k}^2+\frac{P\Theta_k}{n}} +\square\frac{\Upsilon x^2}{\eta^3 n}\enspace,$$
 $$\norm{\bayes-\bayes_k}^2+\frac{P\Theta_k}{n} \leq (1+\eta) \norm{\bayes-\ERM_k}^2+\square\frac{\Upsilon x^2}{\eta^3 n}\enspace.$$
\end{proposition}
\noindent This means that not only in expectation but also with high probability can the term $\norm{\bayes-\ERM_k}^2$ be decomposed in a bias term $\norm{\bayes-\bayes_k}^2$ and a "variance"  term $P\Theta_k/n$. The bias term measures the capacity of the kernel $k$ to approximate $\bayes$ whereas $P\Theta_k/n$ is the price to pay for replacing $\bayes_k$ by its empirical version $\ERM_k$. In this sense, $P\Theta_k/n$ measures the complexity of the kernel $k$  in a way which is completely adapted to our problem of density estimation. Even if it does not seem like a natural measure of complexity at first glance, note that in the previous examples, it is indeed always linked to a natural complexity. When dealing with regular histograms defined on $[0,1]$, $P\Theta_{k_S}$ is the dimension of the considered space $S$, whereas for approximation kernels  $P\Theta_{k_{K,h}}$ is proportional to the inverse of the considered bandwidth $h$.

\subsection{Some general results about the minimal penalty}
In this section, we assume that we are in the asymptotic regime where the number of observations $n\to\infty$. In particular, the asymptotic notations refers to this regime. 

\noindent From now on, the family $\sK=\sK_n$ may depend on $n$ as long as both $\Gamma$ and $\Upsilon$ remain absolute constants that do not depend on it. Indeed, on the previous examples, this seems a reasonable regime. Since $\sK_n$ now depends on $n$, our selected $\kh=\hat{k}_n$ also depends on $n$.

\noindent To prove that the lower bound $\penmin(k)$ is sharp, we need to show that the estimator chosen by minimizing \eqref{eq:CritPen} with a penalty smaller than $\penmin$ does not satisfy an oracle inequality. This is only possible if the $\norm{\bayes-\ERM_k}^2$'s are not of the same order and if they are larger than the remaining term $\square(r+\theta^{-3})\Upsilon x^2/n$. From an asymptotic point of view, we rewrite  this thanks to \propref{ConcRisk} as for all $n\ge 1$, there exist $k_{0,n}$ and $k_{1,n}$ in $\sK_n$ such that
\begin{equation}\label{eq:1stcondMin}
\norm{\bayes-\bayes_{k_{1,n}}}^2\!+\!\frac{P\Theta_{k_{1,n}}}{n} \!\gg\!\norm{\bayes-\bayes_{k_{0,n}}}^2\!+\!\frac{P\Theta_{k_{0,n}}}{n} \! \gg \! \square\paren{\! r+\frac{1}{\theta^3}\!} \frac{\Upsilon x^2}{n} \enspace,
\end{equation}
where $a_n\gg b_n$ means that $b_n/a_n\to_{n\to\infty}0$.
More explicitly, denoting by $\petito(1)$ a sequence only depending on $n$ and tending to 0 as $n$ tends to infinity and whose value may change from line to line, one assumes that there exists $c_\bayes$ and $c_R$ positive constants such that for all $n\ge 1$, there exist $k_{0,n}$ and $k_{1,n}$ in $\sK_n$ such that
\vspace{-0.2cm}
\begin{eqnarray}\label{eq:condMinbis1}
 \norm{\bayes-\bayes_{k_{0,n}}}^2+\frac{P\Theta_{k_{0,n}}}{n}  \leq c_\bayes~ \petito(1) \paren{\norm{\bayes-\bayes_{k_{1,n}}}^2+\frac{P\Theta_{k_{1,n}}}{n}} \\
\frac{(\log(|\sK_n|\vee n))^3}{n} \leq  c_R~ \petito(1)\paren{\norm{\bayes-\bayes_{k_{0,n}}}^2+\frac{P\Theta_{k_{0,n}}}{n}}\enspace.\label{eq:condMinbis2}
\end{eqnarray}
\vspace{-0.3cm}

\noindent We put a log-cube factor in  the remaining term to allow some  choices of $\theta=\theta_n\to_{n\to\infty} 0$ and $r=r_n\to_{n\to\infty} +\infty$.

\noindent But \eqref{eq:condMinbis1} and \eqref{eq:condMinbis2} (or \eqref{eq:1stcondMin}) are not sufficient. Indeed, the following result explains what happens when 
the bias terms are always the leading terms.

\begin{corollary}\label{cor:biaisdomine}
Let $(\sK_n)_{n\geq 1}$ be a sequence of finite collections of kernels $k$ satisfying Assumptions~\eqref{eq:AssOr1}, \eqref{eq:AssOr2}, \eqref{eq:AssOr3}, \eqref{eq:AssOr4} \eqref{eq:AssOr5}, \eqref{eq:AssVar} for a positive constant $\Upsilon$ independent of $n$ and  such that 
\vspace{-0.1cm}
\begin{equation}\label{eq:CondGrosBiais}
 \frac{1}{n} = c_b~\petito(1)\inf_{k\in \sK_n} \frac{\norm{\bayes-\bayes_k}^2}{P\Theta_k} \enspace,
 \end{equation}
 \vspace{-0.3cm}
 
\noindent for some positive constant $c_b$.

\noindent Assume that there exist real numbers of any sign $\delta'\ge \delta$ and a sequence  $(r_n)_{n\geq 1}$ of nonnegative real numbers such that, for all $n\ge 1$, with probability larger than $1-\square / n^2$, for all $k\in \sK_n$,
\vspace{-0.2cm}
\begin{multline*}
\delta\frac{P\Theta_k}n-\square_{\delta,\delta',\Upsilon} \frac{r_n \log(n\vee |\sK_n|)^2}{n}\\
\le \pen(k)-\frac{2P\chi_k-P\Theta_k}{n}\le \delta'\frac{P\Theta_k}n+\square_{\delta,\delta',\Upsilon} \frac{r_n\log(n\vee |\sK_n|)^2}{n}\enspace.
\end{multline*}
 \vspace{-0.3cm}
 
\noindent 
Then, with probability larger than $1-\square/n^2$,
\vspace{-0.3cm}
\begin{multline*}
\norm{\bayes-\ERM_{\kh_n}}^2\le\\
 (1+\square_{\delta,\delta',\Upsilon,c_b}~\petito(1))\inf_{k\in \sK_n}\norm{\bayes-\ERM_k}^2+\square_{\delta,\delta',\Upsilon} \paren{r_n+\log n}\frac{\log(n\vee |\sK_n|)^2}{n}\enspace. 
\end{multline*}
 \vspace{-0.4cm}
\end{corollary}
\noindent  The proof easily follows by taking $\theta=(\log n)^{-1}$ in \eqref{eq:WhateverPenalty}, $\eta=2$ for instance in \propref{ConcRisk} and by  using Assumption~\eqref{eq:CondGrosBiais} and the bounds on $\pen(k)$. This result shows that the estimator $\ERM_{\kh_n}$ satisfies an asymptotically optimal oracle inequality when condition \eqref{eq:CondGrosBiais} holds, whatever the values of $\delta$ and $\delta'$ even when they are negative. This proves that the lower bound $\penmin$ is not sharp in this case. 
%
%

\noindent  Therefore, we have to assume that at least one bias $\norm{\bayes-\ERM_k}^2$ is negligible with respect to $P\Theta_k/n$. Actually, to conclude, we assume that this happens for $k_{1,n}$ in \eqref{eq:1stcondMin}.

\begin{theorem}\label{thm:minpen}
Let  $(\sK_n)_{n\geq 1}$ be a sequence of  finite collections of kernels satisfying Assumptions~\eqref{eq:AssOr1}, \eqref{eq:AssOr2}, \eqref{eq:AssOr3}, \eqref{eq:AssOr4} \eqref{eq:AssOr5}, \eqref{eq:AssVar}, with $\Upsilon$ not depending on $n$. Each $\sK_n$ is also assumed to satisfy \eqref{eq:condMinbis1} and \eqref{eq:condMinbis2} with a kernel $k_{1,n}\in\sK_n$ in \eqref{eq:condMinbis1} such that
\vspace{-0.2cm}
\begin{equation}\label{eq:2ndcondMin}
\norm{\bayes-\bayes_{k_{1,n}}}^2 \leq c~\petito(1) \frac{P\Theta_{k_{1,n}}}{n}\enspace,
\end{equation}
\vspace{-0.3cm}

\noindent for some fixed positive constant $c$.
Suppose that there exist $\delta\geq \delta'>0$ and a sequence  $(r_n)_{n\geq 1}$ of nonnegative real numbers such that $r_n\leq \square\log(|\sK_n|\vee n)$ and such that for all $n\ge 1$, with probability larger than $1-\square  n^{-2}$, for all $k\in \sK_n$,
\vspace{-0.2cm}
\begin{multline}\label{eq:CondPenMin}
\frac{2P\chi_k-P\Theta_k}{n}-\delta\frac{P\Theta_k}n-\square_{\delta,\delta',\Upsilon} \frac{r_n \log(|\sK_n|\vee n)^2}{n}\le \pen(k)\\
\le \frac{2P\chi_k-P\Theta_k}{n}-\delta'\frac{P\Theta_k}n+\square_{\delta,\delta',\Upsilon} \frac{r_n \log(|\sK_n|\vee n)^2}{n}\enspace. 
\end{multline}
\vspace{-0.3cm}

\noindent
Then, with probability larger than $1-\square/n^2$, the following holds
\vspace{-0.1cm}
\begin{equation}\label{eq:Ptheta}
P\Theta_{\kh_n}\ge \paren{\frac{\delta'}{\delta}+\square_{\delta,\delta',\Upsilon,c,c_s,c_R}~\petito(1)}P\Theta_{k_{1,n}} \quad \mbox{and}
\end{equation}
\vspace{-0.5cm}
\begin{multline}
\norm{\bayes-\ERM_{\kh_n}}^2\ge \paren{\frac{\delta'}{\delta}+\square_{\delta,\delta',\Upsilon,c,c_s,c_R}~\petito(1)}\norm{\bayes-\ERM_{k_{1,n}}}^2\\
\gg\norm{\bayes-\ERM_{k_{0,n}}}^2\ge \inf_{k\in \sK_n}\norm{\bayes-\ERM_{k}}^2\enspace.\label{eq:risk}
\end{multline}
\vspace{-0.7cm}
\end{theorem}

\noindent By \eqref{eq:risk}, under the conditions of \thmref{minpen}, the estimator $\ERM_{\kh_n}$ cannot satisfy an oracle inequality, hence, the lower bound $(2P\chi_k-P\Theta_k)/n$ in \thmref{OptimalPenalty} is sharp. This shows that $(2P\chi_k-P\Theta_k)/n$ is a minimal penalty in the sense of \cite{Bir_Mas:2007} for kernel selection. When 
\vspace{-0.1cm}
\[\pen(k)=\frac{2P\chi_k-P\Theta_k}{n}+\kappa \frac{P\Theta_k}n \enspace,\]
\vspace{-0.5cm}

\noindent the complexity $P\Theta_{\kh_n}$ presents a sharp phase transition when $\kappa$ becomes positive. Indeed, when $\kappa<0$ it follows from \eqref{eq:Ptheta} that the complexity $P\Theta_{\kh_n}$ is asymptotically larger than $P\Theta_{k_{1,n}}$. But on the other hand, as a consequence of \thmref{OptimalPenalty}, when $\kappa>0$, this complexity becomes smaller than 
\vspace{-0.1cm}
\begin{multline*}
\square_{\kappa}n\inf_{k\in\sK_n}\paren{\norm{\bayes-\bayes_k}^2+\frac{P\Theta_k}n}\le \square_{\kappa}\paren{n\norm{\bayes-\bayes_{k_{0,n}}}^2+P\Theta_{k_{0,n}}}\\
\ll\square_{\kappa}\paren{n\norm{\bayes-\bayes_{k_{1,n}}}^2+P\Theta_{k_{1,n}}}\le \square_{\kappa} P\Theta_{k_{1,n}}\enspace. 
\end{multline*}
\vspace{-0.3cm}
\subsection{Examples}
\paragraph{Example 1 (continued).}
Let $\sS=\sS_n$ be the collection of spaces of  regular histograms on $[0,1]$ with dimensions $\set{1,\ldots,n}$ and let $\hat{S}=\hat{S}_n$ be the selected space thanks to the penalized criterion. Recall that, for any $S\in\sS_n$, the orthonormal basis is defined by $\vphi_i=\sqrt{D_S}{\bf 1}_{i}$ and $P\Theta_{k_S}=D_S$.
Assume that $\bayes$ is $\alpha$-H\"olderian, with $\alpha\in(0,1]$ with $\alpha$-H\"olderian norm $L$. It is well known (see for instance Section 1.3.3. of \cite{Birge2006}) that the bias is bounded above by
\[\norm{\bayes-\bayes_{k_S}}^2\le \square_L D_S^{-2\alpha}\enspace.\]
\vspace{-0.5cm}

\noindent In particular, if $D_{S_1}=n$,
\[\norm{\bayes-\bayes_{k_{S_1}}}^2\le \square_L n^{-2\alpha}\ll 1= \frac{D_{S_1}}n=\frac{P\Theta_{k_{S_1}}}n\enspace.\]
Thus, \eqref{eq:2ndcondMin} holds for kernel $k_{S_1}$. Moreover, if $D_{S_0}=\lfloor\sqrt{n}\rfloor$, 
\vspace{-0.1cm}
\begin{multline*}
\frac{(\log(n\vee|\sS_n|)^{3}}n\ll\norm{\bayes-\bayes_{k_{S_0}}}^2+\frac{D_{S_0}}n\le \square_L \paren{\frac1{n^\alpha}+\frac1{\sqrt n}}\\
\ll \norm{\bayes-\bayes_{k_{S_1}}}^2+\frac{D_{S_1}}n\enspace. 
\end{multline*}
\vspace{-0.4cm}

\noindent Hence, \eqref{eq:1stcondMin} holds with $k_{0,n}=k_{S_0}$ and $k_{1,n}=k_{S_1}$. Therefore, \thmref{minpen} and \thmref{OptimalPenalty} apply in this example. If $\pen(k_S)= (1-\delta)D_S/n$, the dimension $D_{k_{\widehat{S}_n}}\ge \square_\delta n$ and $\ERM_{k_{\widehat{S}_n}}$ is not consistent and does not satisfy an oracle inequality.  On the other hand, if $\pen(k_S)=(1+\delta)D_S/n$, 
\[D_{{\widehat{S}_n}}\le \square_{L,\delta} \paren{n^{1-\alpha}+\sqrt n}\ll D_{S_1}=n\]
and $\ERM_{k_{\widehat{S}_n}}$ satisfies an oracle inequality which implies that, with probability larger than $1-\square/n^2$,
\[\norm{\bayes-\ERM_{k_{\widehat{S}_n}}}^2\le \square_{\alpha,L,\delta} n^{-2\alpha/(2\alpha+1)}\enspace,\]
by taking $D_S\simeq n^{1/(2\alpha+1)}.$
It achieves the minimax rate of convergence over the class of $\alpha$-H\"olderian functions. 

\noindent  From \thmref{OptimalPenalty}, the penalty $\pen(k_S)=2D_S/n$ provides an estimator $\ERM_{k_{\widehat{S}_n}}$ that achieves an asymptotically optimal oracle inequality. Therefore  the optimal penalty is equal to $2$ times the minimal one. In particular, the slope heuristics of \cite{Bir_Mas:2007} holds in this example, as already noticed in \cite{Le2012}.

\noindent  Finally to illustrate \corref{biaisdomine}, let us take $\bayes(x)=2x$ and  the collection of regular histograms with dimension in $\{1,\ldots,\lfloor n^{\beta}\rfloor\}$, with $\beta<1/3$. Simple calculations show that 
$$\frac{\norm{\bayes-\bayes_{k_S}}^2}{D_S}\geq \square D_S^{-3}\geq \square n^{-3\beta} \gg n^{-1}.$$  
Hence \eqref{eq:CondGrosBiais} applies and the penalized estimator  with penalty $\pen(k_S)\simeq \delta\frac{D_S}{n}$ always satisfies an oracle inequality even if $\delta =0$ or $\delta<0$. This was actually expected since it is likely to choose the largest dimension which is also the oracle choice in this case.

\paragraph{Example 2 (continued).} Let $K$ be a fixed function, let $\sH=\sH_n$ denote the following grid of bandwidths 
\[\sH=\set{\frac{\norm{K}_\infty \norm{K}_1}{i} \quad / \quad i=1,\ldots,n}\enspace\]
and let $\hat{h}=\hat{h}_n$ be the selected bandwidth.
Assume as before that $\bayes$ is a density on $[0,1]$ that belongs to the Nikol'ski class $\mathcal{N}(\alpha,L)$ with $\alpha\in(0,1]$ and $L>0$. By Proposition~1.5 in \cite{Tsy:2009}, if $K$ satisfies $\int \absj{u}^\alpha \absj{K(u)} du<\infty$
\vspace{-0.1cm}
\[\norm{\bayes-\bayes_{k_{K,h}}}^2\le \square_{\alpha,K,L}h^{2\alpha}\enspace.\]
\vspace{-0.3cm}

\noindent In particular, when $h_1=\norm{K}_\infty\norm{K}_1/n$, 
\vspace{-0.1cm}
\[\norm{\bayes-\bayes_{k_{K,h_1}}}^2\le \square_{\alpha,K,L}n^{-2\alpha}\ll \frac{P\Theta_{k_{K,h_1}}}n=\frac{\norm{K}^2}{\norm{K}_\infty\norm{K}_1}\enspace.\]
\vspace{-0.3cm}

\noindent 
On the other hand,  for $h_0=\norm{K}_\infty\norm{K}_1/\PEInf{\sqrt{n}}$, 
\vspace{-0.1cm}
\begin{multline*}
\frac{(\log n\vee|\sH_n|)^2}n\ll\norm{\bayes-\bayes_{k_{K,h_0}}}^2+\frac{P\Theta_{k_{K,h_0}}}n\\
\le \square_{K,\alpha,L}\paren{\frac1{n^\alpha}+\frac1{\sqrt n}}\ll \norm{\bayes-\bayes_{k_{K,h_1}}}^2+\frac{P\Theta_{k_{K,h_1}}}n\enspace. 
\end{multline*}
\vspace{-0.3cm}

\noindent 
Hence,  \eqref{eq:1stcondMin} and \eqref{eq:2ndcondMin} hold with kernels $k_{0,n}=k_{K,h_0}$ and $k_{1,n}=k_{K,h_1}$. Therefore, \thmref{minpen} and \thmref{OptimalPenalty} apply in this example. If for some $\delta>0$ we set $\pen(k_{K,h})= (2K(0)-\norm{K}^2-\delta\norm{K}^2)/(nh)$, then  $\widehat{h}_n\le \square_{\delta,K} n^{-1}$ and $\ERM_{k_{K,\widehat{h}_n}}$ is not consistent and does not satisfy an oracle inequality. On the other hand, if $\pen(k_{K,h})= (2K(0)-\norm{K}^2+\delta\norm{K}^2)/(nh)$, then
\vspace{-0.1cm}
\[\widehat{h}_n\ge \square_{\delta,K,L} \paren{n^{1-\alpha}+\sqrt{n}}^{-1} \gg \square_{\delta,K,L} n^{-1}\enspace,\]
\vspace{-0.3cm}

\noindent 
and $\ERM_{K,k_{\widehat{h}_n}}$ satisfies an oracle inequality which implies that, with probability larger than $1-\square/n^2$,
\vspace{-0.1cm}
\[\norm{\bayes-\ERM_{k_{K,\widehat{h}_n}}}^2\le \square_{\alpha,K,L,\delta} n^{-2\alpha/(2\alpha+1)}\enspace,\]
\vspace{-0.3cm}

\noindent 
for $h=\norm{K}_\infty \norm{K}_1/\PEInf{n^{1/(2\alpha+1)}}\in \sH.$
In particular it achieves the minimax rate of convergence over the class $\mathcal{N}(\alpha,L)$. Finally, if $\pen(k_{K,h})=2K(0)/(nh)$, $\ERM_{k_{K,\widehat{h}_n}}$ achieves an asymptotically optimal oracle inequality, thanks to  \thmref{OptimalPenalty}. 

\noindent The minimal penalty is therefore 
\vspace{-0.2cm}
$$\penmin(k_{K,h})=\frac{2K(0)-\norm{K}^2}{nh}\enspace.$$
\vspace{-0.3cm}

\noindent In this case, the optimal penalty $\pen_{\rm opt}(k_{K,h})=2K(0)/(nh)$ derived from \thmref{OptimalPenalty} is not twice the minimal one, but one still has, if $2K(0)\ne \norm{K}^2$,
\vspace{-0.2cm}
$$\penopt(k_{K,h})=\frac{2K(0)}{2K(0)-\norm{K}^2} \penmin(k_{K,h})\enspace,$$
\vspace{-0.2cm}

\noindent
even if they can  be of opposite sign depending on $K$. This type of nontrivial relationship between optimal and minimal penalty has already been underlined in \cite{Arl_Bac:2009} in regression framework for selecting linear estimators.

\noindent Note that if one allows two kernel functions $K_1$ and $K_2$ in the family of kernels such that $2K_1(0)\ne \norm{K_1}^2$, $2K_2(0)\ne \norm{K_2}^2$ and
\[\frac{2K_1(0)}{2K_1(0)-\norm{K_1}^2}\ne \frac{2K_2(0)}{2K_2(0)-\norm{K_2}^2}\enspace,\]
then there is no absolute constant multiplicative factor linking the minimal penalty and the optimal one. 
\vspace{-0.5cm}
\section{Small simulation study\label{sec:simu}}
In this section we illustrate on simulated data \thmref{OptimalPenalty} and \thmref{minpen}.  We focus on approximation kernels only, since projection kernels have been already discussed in \cite{Le2012}. 

\noindent We observe an $n=100$ i.i.d. sample of standard gaussian distribution. For a fixed parameter $a\ge 0$ we consider the family of kernels 
\vspace{-0.1cm}
$$k_{K_a,h}(x,y)=\frac1hK_a\paren{\frac{x-y}h}\qquad \mbox{with}\quad h\in \sH=\set{\frac1{2i},~i=1,\ldots,50}\enspace,$$
\vspace{-0.2cm}

\noindent 
where for $\displaystyle x\in \R,\quad
K_a(x)=\frac{1}{2\sqrt{2\pi}} \paren{e^{-\frac{(x-a)^2}{2}}+e^{-\frac{(x+a)^2}{2}}}\enspace.$\\
In particular the kernel estimator with $a=0$ is the classical Gaussian kernel estimator. Moreover
\vspace{-0.1cm}
$$K_a(0)=\frac1{\sqrt{2\pi}}\exp\paren{-\frac{a^2}2} \quad \mbox{and}\quad \norm{K_a}^2=\frac{1+e^{-a^2}}{4\sqrt{\pi}}\enspace.$$
\vspace{-0.3cm}

\begin{figure}[h!]
\begin{tabular}{c}
~\vspace{-1cm}\\
\hspace{-1cm}\includegraphics[width=10cm,height=4cm]{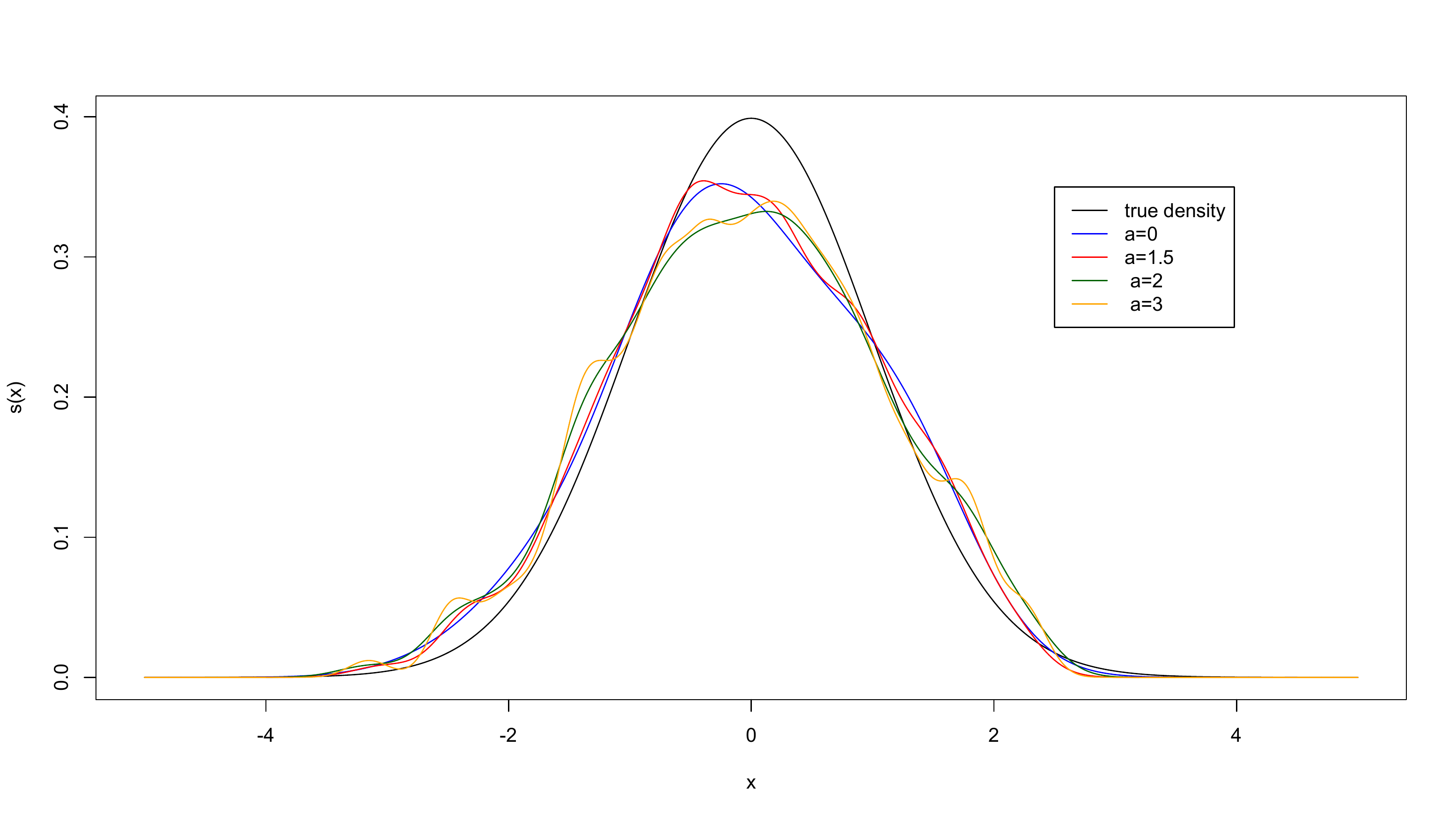}
\end{tabular}
\caption{\label{fig:optimal} {\small Selected approximation kernel estimators when the penalty is the optimal one i.e. $\frac{2K_a(0)}{nh}$}. }
\end{figure}

\begin{figure}[h!]
\begin{tabular}{c}
\hspace{-1cm}\includegraphics[width=13cm,height=4cm]{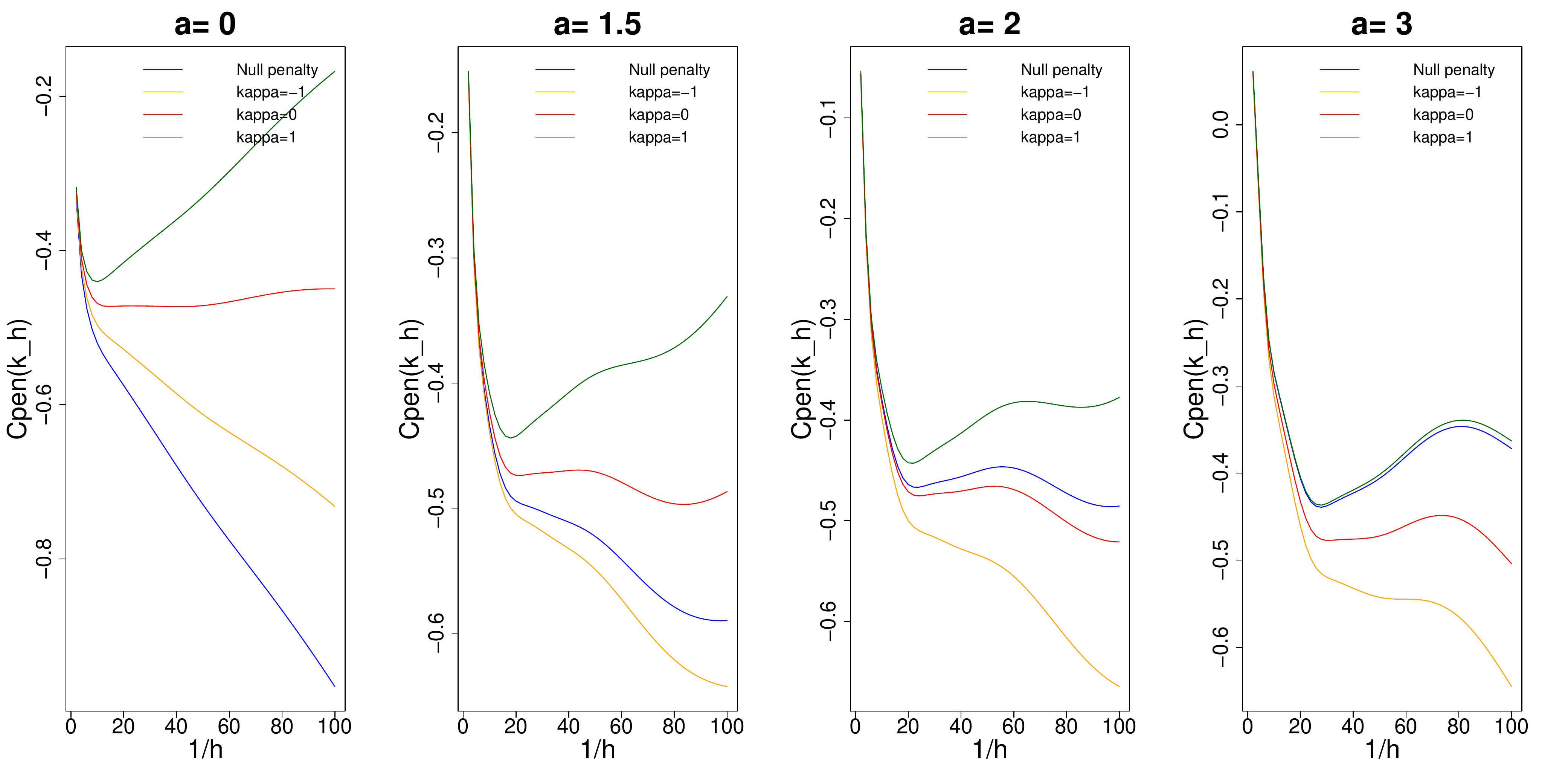}
\end{tabular}
\caption{\label{fig:contraste} {\small Behavior of $P_n\gamma(\ERM_{k_{K_a,h}})$ (blue line) and $\sC_{\pen}\paren{k_{K_a,h}}$ as a function of $1/h$, which is proportional to the complexity $P\Theta_{k_{K_a,h}}$.}\vspace{-0.5cm} }
\end{figure}

\noindent 
Thus, depending on the value of $a$, the minimal penalty $(2K_a(0)-\norm{K_a}^2)/(nh)$ may be negative. We study the behavior of the penalized criterion 
\vspace{-0.1cm}
\[\sC_{\pen}\paren{k_{K_a,h}}=  P_n\gamma(\ERM_{k_{K_a,h}})+\pen(k_{K_a,h})\]
\vspace{-0.3cm}

\noindent 
with penalties of the form
\begin{equation}\label{eq:penkappa}
\pen\paren{k_{K_a,h}}=\frac{2K_a(0)-\norm{K_a}^2}{nh}+\kappa \frac{\norm{K_a}^2}{nh}\enspace, 
\end{equation}
for different values of $\kappa$ ($\kappa=-1,0,1$) and $a$ ($a=0,1.5,2,3$). On \figref{optimal} are represented the selected estimates by the optimal penalty $2K_a(0)/(nh)$ for the different values of $a$ and on \figref{contraste} one sees the evolution of the different penalized criteria as a function of $1/h$.
The contrast curves for $a=0$ are classical on  \figref{contraste}. Without penalization, the criterion decreases and leads to the selection of the smallest bandwidth. At the minimal penalty, the curve is flat and at the optimal penalty one selects a meaningful bandwidth as shown on \figref{optimal}.

\noindent When $a>0$, despite the choice of those unusual kernels, the reconstructions on \figref{optimal} for the optimal penalty are also meaningful. However when $a=2$ or $a=3$, the criterion with minimal penalty is smaller than the unpenalized criterion, meaning that minimizing the latter criterion leads by \thmref{OptimalPenalty} to an oracle inequality. In our simulation, when $a=3$, the curves for the optimal criterion and the unpenalized one are so close that the same estimator is selected by both methods.

\begin{figure}[h!]
\begin{tabular}{c}
\hspace{-1cm}\includegraphics[width=13cm,height=4cm]{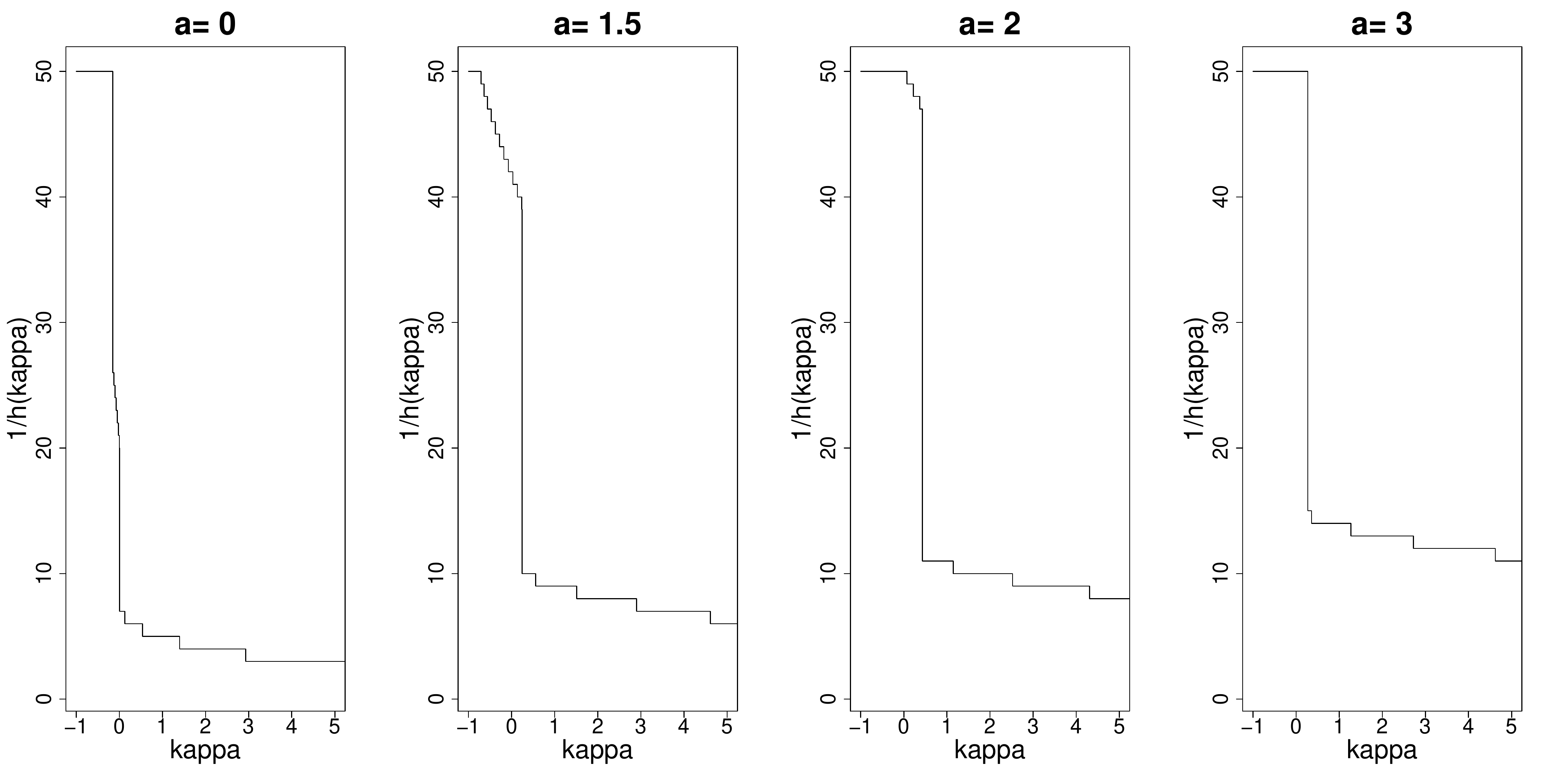}
\end{tabular}
\caption{\label{fig:saut}  {\small Behavior of  $1/\hat{h}$, which is proportional to the complexity $P\Theta_{k_{K_a,h}}$, for the estimator selected by the criterion whose penalty is given by \eqref{eq:penkappa}, as a function of $\kappa$.} ~\vspace{-0.3cm}}
\end{figure}
\noindent Finally \figref{saut} shows that there is indeed in all cases a sharp phase transition around $\kappa=0$ i.e. at the minimal penalty for the complexity of the selected estimate.

\section{Proofs \label{sec:proofs}}
\subsection{Proof of \thmref{OptimalPenalty}}
The starting point to prove the oracle inequality is to notice that any minimizer $\kh$ of $\sC_{\pen}$ satisfies
\[\norm{\bayes-\ERM_{\kh}}^2\le \norm{\bayes-\ERM_k}^2+\paren{\pen(k)-\penid(k)}-\paren{\pen\paren{\kh}-\penid\paren{\kh}}\enspace.\]
Using the expression of the ideal penalty \eqref{eq:penid} we find
\vspace{-0.2cm}
\begin{align}\label{eq:DecRisk}
\notag \norm{\bayes-\ERM_{\kh}}^2&\le \norm{\bayes-\ERM_k}^2+\paren{\pen(k)-2\frac{P\chi_k}n}-\paren{\pen\paren{\kh}-2\frac{P\chi_{\kh}}n}\\
\notag &+2\frac{P(\bayes_k-\bayes_{\kh})}n+2\paren{1-\frac2n}(P_n-P)(\bayes_{\kh}-\bayes_k)\\
 &+2\frac{(P_n-P)(\chi_{\kh}-\chi_k)}n+2\frac{U_{\kh}-U_k}{n^2}\enspace.
\end{align}
\vspace{-0.3cm}

\noindent By \propref{ConcRemainder} (see the appendix), for all $x>1$, for all $\theta$ in $(0,1)$, with probability larger than $1-(7.4|\sK|+2|\sK|^2)e^{-x}$,
\vspace{-0.1cm}
\begin{align*}
 \norm{\bayes-\ERM_{\kh}}^2&\le \norm{\bayes-\ERM_k}^2+\paren{\pen(k)-2\frac{P\chi_k}n}-\paren{\pen\paren{\kh}-2\frac{P\chi_{\kh}}n}\\
 & +\theta\norm{\bayes-\bayes_\kh}^2+\theta\norm{\bayes-\bayes_k}^2+ \square \frac{\Upsilon}{\theta n}\\
& + \paren{1-\frac2n} \theta \norm{\bayes-\bayes_\kh}^2 + \paren{1-\frac2n} \theta \norm{\bayes-\bayes_k}^2+  \square \frac{\Upsilon x^2}{\theta n} \\
 & +\theta\frac{P\Theta_k}{n} + \theta\frac{P\Theta_\kh}{n} + \square \frac{\Upsilon x}{\theta n}  + \theta\frac{P\Theta_k}{n} + \theta\frac{P\Theta_\kh}{n} + \square \frac{\Upsilon x^2}{\theta n}
\end{align*}
\vspace{-0.3cm}

\noindent 
Hence
\vspace{-0.3cm}
\begin{align*}
 \norm{\bayes-\ERM_{\kh}}^2
&\leq  \norm{\bayes-\ERM_k}^2+\paren{\pen(k)-2\frac{P\chi_k}n}-\paren{\pen\paren{\kh}-2\frac{P\chi_{\kh}}n}\\
&+2\theta\croch{\norm{\bayes-\bayes_\kh}^2+\frac{P\Theta_\kh}{n}} +2\theta\croch{\norm{\bayes-\bayes_k}^2+\frac{P\Theta_k}{n}} +  \square \frac{\Upsilon x^2}{\theta n}\enspace.
\end{align*}
\vspace{-0.3cm}

\noindent 
This bound holds using \eqref{eq:AssOr1}, \eqref{eq:AssOr2} and \eqref{eq:AssOr3} only. Now by \propref{ConcRisk} applied with $\eta=1$, we have for all $x>1$, for all $\theta\in (0,1)$, with probability larger than $1-(16.8|\sK|+2|\sK|^2) e^{-x}$,
\vspace{-0.1cm}
\begin{align*}
 \norm{\bayes-\ERM_{\kh}}^2&\le \norm{\bayes-\ERM_k}^2+\paren{\pen(k)-2\frac{P\chi_k}n}-\paren{\pen\paren{\kh}-2\frac{P\chi_{\kh}}n}\\
&+ 4 \theta \norm{\bayes-\ERM_\kh}^2 + 4 \theta \norm{\bayes-\ERM_k}^2 + \square \frac{\Upsilon x^2}{\theta n}\enspace.
\end{align*}
\vspace{-0.3cm}

\noindent 
This gives the first part of the theorem. 

\medskip

\noindent For the second part, by the condition~\eqref{eq:CondPen} on the penalty, we find for all $x>1$, for all $\theta$ in $(0,1)$, with probability larger than $1-(C+16.8|\sK|+2|\sK|^2)e^{-x}$,
\begin{multline*}
 (1-4\theta) \norm{\bayes-\ERM_{\kh}}^2\le \\
 (1+4\theta)\norm{\bayes-\ERM_k}^2+ (\delta'-1)_+ \frac{P\Theta_k}{n} + (1-\delta)_+ \frac{P\Theta_\kh}{n} + \square \paren{r+\frac1{\theta}} \frac{\Upsilon x^2}{n}\enspace. 
\end{multline*}
By  \propref{ConcRisk} applied with $\eta=\theta$, we have with probability larger than $1-(C+26.2|\sK|+2|\sK|^2)e^{-x}$,
\vspace{-0.2cm}
\begin{multline*}
 (1-4\theta) \norm{\bayes-\ERM_{\kh}}^2\le (1+4\theta)\norm{\bayes-\ERM_k}^2+ (\delta'-1)_+ (1+\theta)  \norm{\bayes-\ERM_k}^2\\
  + (1-\delta)_+ (1+\theta)  \norm{\bayes-\ERM_{\kh}}^2 + \square \paren{r+\frac1{\theta^3}} \frac{\Upsilon x^2}{n}\enspace,
 \end{multline*}
 \vspace{-0.6cm}
 
\noindent that is
\vspace{-0.1cm}
\begin{multline*}
 \paren{(\delta\wedge1)-\theta(4+(1-\delta)_+)} \norm{\bayes-\ERM_{\kh}}^2\\
 \le \paren{(\delta'\vee 1)+\theta(4+(\delta'-1)_+)} \norm{\bayes-\ERM_k}^2+ \square \paren{r+\frac1{\theta^3}} \frac{\Upsilon x^2}{n}\enspace. 
  \end{multline*}
Hence, because $1\leq [(\delta'\vee 1)+(4+(\delta'-1)_+)\theta]\leq  (\delta'\vee 1)+(4+\delta')\theta$, we obtain the desired result.
\subsection{Proof of \propref{ConcRisk}}
%

First, let us denote for all $x\in \Xbb$ 
\vspace{-0.1cm}
\[F_{A,k}(x)\egaldef \Ex{A_k(X,x)}, \qquad \zeta_k(x)\egaldef\int\paren{k(y,x)-\bayes_k(y)}^2d\mu(y)\enspace,\]
\vspace{-0.3cm}

\noindent
and 
\vspace{-0.1cm}
 \[U_{A,k} \egaldef \sum_{i\ne j=1}^n \paren{A_k(X_i,X_j)-F_{A,k}(X_i)-F_{A,k}(X_j)+\Ex{A_k(X,Y)}}\enspace.\]
 \vspace{-0.3cm}

\noindent
Some easy computations then provide the following useful equality 
 \[\norm{\bayes_k-\ERM_k}^2 = \frac1nP_n\zeta_k+\frac1{n^2}U_{A,k}\enspace.\]

\noindent We need only treat the terms on the right-hand side, thanks to the probability tools of \secref{tools}. Applying \propref{Bern}, we get, for any $x\ge 1$, with probability larger than $1-2\absj{\sK}e^{-x}$,
\[\absj{(P_n-P)\zeta_k}\le\sqrt{\frac{2x}n P\zeta_k^2}+\frac{\norm{\zeta_k}_\infty x}{3n}\enspace.\]
One can then check the following link between $\zeta_k$ and $\Theta_k$
\[P\zeta_k=\int\paren{k(y,x)-\bayes_k(x)}^2\bayes(y)d\mu(x)d\mu(y)=P\Theta_k-\norm{\bayes_k}^2 \enspace.\]
Next, by \eqref{eq:CdtForGamma} and \eqref{eq:AssOr1}
\begin{align*}
\norm{\zeta_k}_\infty&=\sup_{y\in\Xbb}\int\paren{k(y,x)-\Ex{k(X,x)}}^2d\mu(x)\\
&\le 4\sup_{y\in\Xbb}\int k(y,x)^2d\mu(x)\le 4\Upsilon n\enspace. 
\end{align*}
In particular, since $\zeta_k\ge 0$,
\[ P\zeta_k^2\le \norm{\zeta_k}_\infty P\zeta_k\le 4\Upsilon nP\Theta_k\enspace.\]
It follows from these computations and from \eqref{eq:AssOr1} that there exists an absolute constant $\square$ such that, for any $x\ge 1$, with probability larger than $1-2\absj{\sK}e^{-x}$, for any $\theta\in(0,1)$,
\vspace{-0.1cm}
\[\absj{P_n\zeta_k-P\Theta_k}\le \theta P\Theta_k+\square \frac{\Upsilon x}{\theta}\enspace.\]
\vspace{-0.3cm}

\noindent We now need to control the term $U_{A,k}$. From \propref{ConcUstat}, for any $x\ge 1$, with probability larger than $1-5.4\absj{\sK}e^{-x}$, 
\[\frac{\absj{U_{A,k}}}{n^2}\le \frac{\square}{n^2}\paren{C\sqrt{x}+Dx+Bx^{3/2}+Ax^2}\enspace.\]
By \eqref{eq:CdtForGamma}, \eqref{eq:AssOr1} and Cauchy-Schwarz inequality,
\[A=4\sup_{(x,y)\in\Xbb^2}\int k(x,z)k(y,z)d\mu(z)\le 4\sup_{x\in\Xbb}\int k(x,z)^2d\mu(z)\le 4\Upsilon n\enspace.\]
In addition, by \eqref{eq:AssOr5},
$B^2\le 16\sup_{x\in\Xbb}\Ex{A_k(X,x)^2}\le 16\Upsilon n\enspace.$\\
Moreover, applying the Assumption~\eqref{eq:AssOr4},
\[C^2\le \sum_{i\ne j=1}^n\Ex{A_k(X_i,X_j)^2}=n^2\Ex{A_k(X,Y)^2}\le n^2 \Upsilon P\Theta_k\enspace.\]
Finally, applying the Cauchy-Schwarz inequality and proceeding as for $C^2$, the quantity used to define $D$ can be bounded above as follows:
\[\Ex{\sum_{i=1}^{n-1}\sum_{j=i+1}^na_i(X_i)b_j(X_j)A_k(X_i,X_j)}\le n\sqrt{\Ex{A_k(X,Y)^2}}\le n\sqrt{\Upsilon P\Theta_k}\enspace.\]
Hence for any $x\ge 1$, with probability larger than $1-5.4\absj{\sK}e^{-x}$, 
$$\mbox{for any } \theta\in(0,1),\quad\frac{\absj{U_{A,k}}}{n^2}\le \theta\frac{P\Theta_k}n+\square \frac{\Upsilon x^2}{\theta n} \enspace.$$
Therefore, for all $\theta\in (0,1)$,
\vspace{-0.2cm}
$$\absj{\norm{\ERM_k-\bayes_k}^2-\frac{P\Theta_k}{n}}\le 2\theta\frac{P\Theta_k}{n}+ \square \frac{\Upsilon x^2}{\theta n} \enspace,$$
\vspace{-0.4cm}

\noindent and the first part  of the result follows by choosing $\theta=\eta/2$.
Concerning the two remaining inequalities appearing in the proposition, we begin by developing the loss. For all $\kK$
$$\norm{\ERM_k-\bayes}^2=\norm{\ERM_k-\bayes_k}^2+\norm{\bayes_k-\bayes}^2+ 2 \psh{\ERM_k-\bayes_k}{\bayes_k-\bayes}\enspace.$$
Then, for all $x\in \Xbb$
\vspace{-0.2cm}
\begin{align*}
 F_{A,k}(x)-\bayes_k(x)&=\int \bayes(y)\int k(x,z)k(z,y) d\mu(z)d\mu(y)-\int \bayes(z)k(z,x)d\mu(z)\\
 &=\int \paren{\int \bayes(y)k(z,y)d\mu(y)-\bayes(z)}k(x,z)d\mu(z)\\
 &=\int\paren{\bayes_k(z)-\bayes(z)}k(z,x)d\mu(z)\enspace.
\end{align*}
\vspace{-0.4cm}

\noindent Moreover, since $PF_{A,k}=\norm{\bayes_k}^2$, we find 
\vspace{-0.2cm}
\begin{align*}
\psh{\ERM_k-\bayes_k}{\bayes_k-\bayes} &= \int \paren{\ERM_k(x)\paren{\bayes_k(x)-\bayes(x)}}d\mu(x)+\Ex{\bayes_k(X)} -\norm{\bayes_k}^2\\
&=\frac{1}{n}\sum_{i=1}^n \int \paren{k(x,X_i)\paren{\bayes_k(x)-\bayes(x)}}d\mu(x)+P(\bayes_k-F_{A,k})\\
&=\frac{1}{n}\sum_{i=1}^n \paren{F_{A,k}(X_i)-\bayes_k(X_i)}+P(\bayes_k-F_{A,k})\\
&= (P_n-P)(F_{A,k}-\bayes_k)\enspace.
\end{align*}
\vspace{-0.4cm}

\noindent 
This expression motivates us to apply again \propref{Bern} to this term. We find by \eqref{eq:CdtForGamma}, \eqref{eq:AssOr1} and  Cauchy-Schwarz inequality
\vspace{-0.2cm}
\begin{align*}
\sup_{x\in\Xbb}\absj{F_{A,k}(x)-\bayes_k(x)}&\le \norm{\bayes-\bayes_k}\sup_{x\in\Xbb}\int \frac{\absj{\bayes(z)-\bayes_k(z)}}{\norm{\bayes-\bayes_k}}k(x,z)d\mu(z)\\
&\le \norm{\bayes-\bayes_k}\sqrt{\sup_{x\in\Xbb}\int k(x,z)^2d\mu(z)}\le \norm{\bayes-\bayes_k}\sqrt{\Upsilon n}\enspace.
\end{align*}
\vspace{-0.5cm}

\noindent 
Moreover,
\vspace{-0.4cm}
\begin{align*}
 P\paren{F_{A,k}-\bayes_k}^2&\le \norm{\bayes-\bayes_k}^2P\paren{\int  \frac{\absj{\bayes(z)-\bayes_k(z)}}{\norm{\bayes-\bayes_k}}k(.,z)d\mu(z)}^2\\
 &\le \norm{\bayes-\bayes_k}^2v_k^2\enspace.
\end{align*}
\vspace{-0.5cm}

\noindent 
Thus by \eqref{eq:AssVar}, for any $\theta,u>0$,
\vspace{-0.1cm}
\begin{align*}
 \sqrt{\frac{2P\paren{F_{A,k}-\bayes_k}^2x}n}\le \theta \norm{\bayes-\bayes_k}^2+\frac{\paren{\Upsilon\vee \sqrt{\Upsilon P\Theta_k}} x}{2\theta n}\\
 \le \theta \norm{\bayes-\bayes_k}^2+\frac{\Upsilon x}{\theta n}\vee \paren{\frac u\theta\frac{ P\Theta_k}{n}+ \frac{\Upsilon x^2}{16\theta u n}}\enspace.
\end{align*}
\vspace{-0.3cm}

\noindent 
Hence, for any $\theta\in(0,1)$ and $x\ge 1$, taking $u=\theta^2$
\begin{align*}
 \sqrt{\frac{2P\paren{F_{A,k}-\bayes_k}^2x}n}\le \theta \paren{\norm{\bayes-\bayes_k}^2+\frac{ P\Theta_k}{n}}+ \square\frac{\Upsilon x^2}{\theta^3 n}\enspace.
\end{align*}
By \propref{Bern}, for all $\theta$ in $(0,1)$ , for all $x>0$ with probability larger than $1-2|\sK| e^{-x}$,
\vspace{-0.3cm}
\begin{align*}
 2\absj{ \psh{\ERM_k-\bayes_k}{\bayes_k-\bayes}} & \leq 2\sqrt{\frac{2P\paren{F_{A,k}-\bayes_k}^2x}n} + 2\norm{\bayes-\bayes_k}\sqrt{\Upsilon n}\frac{x}{3n}\\
& \leq  3\theta \paren{\norm{\bayes-\bayes_k}^2+\frac{ P\Theta_k}{n}}+ \square\frac{\Upsilon x^2}{\theta^3 n}\enspace.
 \end{align*}
 \vspace{-0.4cm}

\noindent Putting together all of the above, one concludes that for all $\theta$ in $(0,1)$, for all $x>1$, with probability larger than $1-9.4|\sK|e^{-x}$
\[
\norm{\ERM_k-\bayes}^2 -\norm{\bayes_k-\bayes}^2\le  3\theta \norm{\bayes-\bayes_k}^2+(1+4\theta)\frac{ P\Theta_k}{n} + \square\frac{\Upsilon x^2}{\theta^3 n} 
\]
and
\[
\norm{\ERM_k-\bayes}^2 -\norm{\bayes_k-\bayes}^2 \ge -3 \theta \paren{\norm{\bayes-\bayes_k}^2+\frac{ P\Theta_k}{n}}+(1-\theta)\frac{ P\Theta_k}{n} - \square\frac{\Upsilon x^2}{\theta^3 n} \enspace.
\]
Choosing, $\theta=\eta/4$ leads to the second part of the result.


\subsection{Proof of \thmref{minpen}}\label{sec:ProofThmMinPen}

It follows from \eqref{eq:WhateverPenalty} (applied with $\theta=\square(\log n)^{-1}$ and $x=\square \log(n\vee |\sK_n|)$) and Assumption~\eqref{eq:CondPenMin} that with probability larger than $1-\square n^{-2}$ we have for any $\kK$ and any $n\geq 2$
\vspace{-0.1cm}
\begin{multline}\label{eq:Cont1}
\norm{\ERM_{\kh_n}-\bayes}^2\leq \paren{1+\frac{\square}{\log n}}\norm{\ERM_k-\bayes}^2-(1+\delta')\paren{1+\frac{\square}{\log n}}\frac{P\Theta_k}n\\
+(1+\delta)\paren{1+\frac{\square}{\log n}}\frac{P\Theta_{\kh_n}}n+\square_{\delta,\delta',\Upsilon}\frac{ \log(|\sK_n|\vee n)^{3}}{n}\enspace.
\end{multline}
\vspace{-0.5cm}

\noindent Applying this inequality with $k=k_{1,n}$ and using \propref{ConcRisk} with $\eta=\square (\log n)^{-1/3}$ and $x=\square\log(|\sK_n|\vee n)$ as a lower bound for $\norm{\ERM_{\kh_n}-\bayes}^2$ and as an upper bound for $\norm{\ERM_{k_{1,n}}-\bayes}^2$, we obtain asymptotically that with probability larger than $1-\square n^{-2}$,
\vspace{-0.1cm}
\begin{multline*}
\notag -\delta (1+\square_\delta ~\petito(1))\frac{P\Theta_{\kh_n}}n\leq \paren{1+\petito(1)}\norm{\bayes_{k_{1,n}}-\bayes}^2-\delta'(1+\square_{\delta'}~\petito(1))\frac{P\Theta_{k_{1,n}}}n\\
+\square_{\delta,\delta',\Upsilon}\frac{\log(|\sK_n|\vee n)^{3}}{n}\enspace.
\end{multline*}
\vspace{-0.4cm}

\noindent By Assumption~\eqref{eq:2ndcondMin}, $\norm{\bayes_{k_{1,n}}-\bayes}^2 \leq  c~ \petito(1) \frac{P\Theta_{k_{1,n}}}n$  and by \eqref{eq:condMinbis1},
\vspace{-0.1cm}
$$\frac{  \paren{\log(|\sK_n|\vee n)}^{3}}{n} \leq c_R c_s ~\petito(1) \frac{P\Theta_{k_{1,n}}}n\enspace.$$ 
\vspace{-0.4cm}

\noindent
This gives \eqref{eq:Ptheta}.
In addition, starting with the event where \eqref{eq:Cont1} holds  and using  \propref{ConcRisk}, we also have with probability larger than $1-\square n^{-2}$,
\vspace{-0.1cm} 
\begin{multline*}
\notag \norm{\ERM_{\kh_n}-\bayes}^2\leq \paren{1+\frac{\square}{\log n}}\norm{\ERM_{k_{1,n}}-\bayes}^2-(1+\delta')\frac{P\Theta_{k_{1,n}}}n\\
+(1+\delta)\paren{1+\petito(1)} \norm{\ERM_{\kh_n}-\bayes}^2+\square_{\delta,\delta',\Upsilon}\frac{ \log(|\sK_n|\vee n)^{3}}{n}\enspace.
\end{multline*}
\vspace{-0.4cm}

\noindent Since $\norm{\ERM_{k_{1,n}}-\bayes}^2 \simeq \frac{P\Theta_{k_{1,n}}}n$, this leads to
\vspace{-0.1cm}
\begin{multline*}
(-\delta+\square_\delta~\petito(1))\norm{\ERM_{\kh}-\bayes}^2\leq \\
-(\delta'+\square_{\delta',c}~\petito(1))\norm{\ERM_{k_{1,n}}-\bayes}^2+\square_{\delta,\delta',\Upsilon}~\frac{ \log(|\sK_n|\vee n)^{3}}{n}\enspace. 
\end{multline*}
\vspace{-0.3cm}

\noindent This leads to \eqref{eq:risk} by \eqref{eq:1stcondMin}.

\vspace{-0.3cm}

\appendix
\section{Proofs for the examples}\label{sec:ProofsKernels}
\subsection{Computation of the constant $\Gamma$ for the three examples}

We have to show for each family $\set{k}_{\kK}$ (see \eqref{eq:DefAk} and \eqref{eq:CdtForGamma}) that there exists a constant $\Gamma\ge 1$ such that for all $\kK$
\[\sup_{x\in\Xbb}~\absj{\Theta_k(x)}\le\Gamma n, \quad \text{and}\quad \sup_{(x,y)\in\Xbb^2}\absj{k(x,y)}\le\Gamma n\enspace.
 \]

\paragraph{Example 1: Projection kernels.}
First, notice that from Cauchy-Schwarz inequality we have for all $(x,y)\in\Xbb^2$ $\absj{k_S(x,y)}\le\sqrt{\chi_{k_S}(x)\chi_{k_S}(y)}$ and by orthonormality, for any $(x,x')\in\Xbb^2$,
\vspace{-0.1cm}
\begin{equation*}
A_{k_S}(x,x')=\sum_{(i,j)\in\sI^2_S}\vphi_i(x)\vphi_j(x')\int_\Xbb \vphi_i(y)\vphi_j(y)d\mu(y)=k_S(x,x')\enspace. 
\end{equation*}
\vspace{-0.3cm}

\noindent In particular, for any $x\in\Xbb$, $\Theta_{k_S}(x)=\chi_{k_S}(x)$. Hence, projection kernels satisfy \eqref{eq:CdtForGamma} for $\Gamma=1\vee n^{-1}\sup_{S\in\sS}\norm{\chi_{k_S}}_\infty$. We conclude by writing
\vspace{-0.2cm}
\[\norm{\chi_{k_S}}_\infty =\sup_{x\in \Xbb} \sum_{i\in\sI_S} \vphi_i(x)^2=\sup_{\substack{(a_i)_{i\in \sI}\telque\\ \sum_{i\in\sI_S} a_i^2=1}} \sup_{x\in \Xbb} \paren{\sum_{i\in\sI_S} a_i\vphi_i(x)}^2\enspace.
\]
For  $f\in S$ we have 
$\norm{f}^2=\sum_{i\in\sI}\psh{f}{\vphi_i}^2$. Hence with $a_i=\psh{f}{\vphi_i}$,
\vspace{-0.1cm}
$$ \norm{\chi_{k_S}}_\infty = \sup_{f\in S, \norm{f}=1} \norm{f}_\infty^2\enspace.$$

\paragraph{Example 2: Approximation kernels.}
$\!\!$First, 
$\sup_{(x,y)\in\Xbb^2}\absj{k_{K,h}(x,y)}\le \norm{K}_\infty/h.$
Second, since  $K\in L^1$
 \vspace{-0.1cm}
\[\Theta_{k_{K,h}}(x)=\frac1{h^2}\int_\Xbb K\paren{\frac{x-y}{h}}^2dy=\frac{\norm{K}^2}{h}\le \frac{\norm{K}_\infty\norm{K}_1}{h} \enspace.\]
\vspace{-0.3cm}

\noindent Now  $K\in L^1$ and $\int K(u) du=1$ implies $\norm{K}_1\ge 1$, hence \eqref{eq:CdtForGamma} holds with $\Gamma=1$ if one assumes that $h\ge \norm{K}_\infty\norm{K}_1/n$.

\paragraph{Example 3: Weighted projection kernels.}
For all $x\in \Xbb$
\vspace{-0.1cm}
\[
\Theta_{k_w}(x)=\sum_{i,j=1}^pw_i\vphi_i(x)w_j\vphi_j(x)\int_\Xbb\vphi_i(y)\vphi_j(y)d\mu(y)= \sum_{i=1}^pw_i^2\vphi_i(x)^2\enspace.
\]
\vspace{-0.3cm}

\noindent From Cauchy-Schwarz inequality, for any $(x,y)\in\Xbb^2$, 
$$\absj{k_w(x,y)}\le \sqrt{\Theta_{k_w}(x)\Theta_{k_w}(y)}\enspace.$$
We thus find that $k_w$ verifies \eqref{eq:CdtForGamma} with $\Gamma \ge 1\vee n^{-1} \sup_{w\in \sW} \norm{\Theta_{k_w}}_\infty$. Since $w_i\le 1$ we find the announced result which is independent of $\sW$.


\subsection{Proof of \propref{AssProjKern}}
Since $\norm{\bayes_{k_S}}^2\le \norm{\bayes}^2\le \norm{\bayes}_\infty$, we find that \eqref{eq:AssOr1} only requires $\Upsilon\ge \Gamma(1+\norm{\bayes}_\infty)$. Assumption~\eqref{eq:AssOr2} holds: this follows from $\Upsilon\ge \Gamma$ and 
\[\Ex{\chi_{k_S}(X)^2}\le \norm{\chi_{k_S}}_\infty P\chi_{k_S} \le \Gamma n P\Theta_{k_S} \enspace.\] 
Now for proving Assumption~\eqref{eq:AssOr4}, we write
\begin{align*}
 \Ex{A_{k_S}(X,Y)^2}&=\Ex{k_S(X,Y)^2}=\int_\Xbb\Ex{k_S(X,x)^2}\bayes(x)d\mu(x)\\
 &\le \norm{\bayes}_\infty\sum_{(i,j)\in\sI^2_S}\Ex{\vphi_i(X)\vphi_j(X)}\int_{\Xbb}\vphi_i(x)\vphi_j(x)d\mu(x)\\
 &=\norm{\bayes}_\infty P\Theta_{k_S}\le \Upsilon  P\Theta_{k_S} \enspace.
\end{align*}
In the same way, Assumption~\eqref{eq:AssOr5} follows from $\norm{s}_\infty \Gamma\le \Upsilon$. Suppose \eqref{eq:Nested} holds with $S=S+S'$ so that the basis $(\vphi_i)_{i\in\sI}$ of $S'$ is included in the one $(\vphi_i)_{i\in\sJ}$ of $S$. Since $\norm{\chi_{k_S}}_\infty\leq \Gamma n$ we have
\begin{align*}
\bayes_{k_S}(x)-\bayes_{k_{S'}}(x)&=\sum_{j\in \sJ\setminus \sI}\paren{P\vphi_j}\vphi_j(x)\leq \sqrt{\sum_{j\in \sJ\setminus \sI}\paren{P\vphi_j}^2\sum_{j\in \sJ\setminus \sI}\vphi_j(x)^2}\\
&\leq \norm{\bayes_{k_S}-\bayes_{k_{S'}}}\norm{\chi_{k_S}}^{1/2}_{\infty} \leq \norm{\bayes_{k_S}-\bayes_{k_{S'}}}\sqrt{\Gamma n}\enspace.
\end{align*}
Hence, \eqref{eq:AssOr3} holds in this case. Assuming \eqref{eq:UnifSupBound} implies that \eqref{eq:AssOr3} holds since
\[\norm{\bayes_{k_S}-\bayes_{k_{S'}}}_\infty\leq \norm{\bayes_{k_S}}_\infty+\norm{\bayes_{k_{S'}}}_\infty\leq \Upsilon \enspace.\]
Finally for \eqref{eq:AssVar}, for any $a\in L^2$, 
\[\int_{\Xbb}a(x)k_S(x,y)d\mu(x)=\sum_{i\in \sI}\psh{a}{\vphi_i}\vphi_i(y)=\Pi_{S}(a)\enspace.\]
is the orthogonal projection of $a$ onto $S$. Therefore, $\Boule_{k_S}$ is the unit ball in $S$ for the $L^2$-norm and, for any $t\in \Boule_{k_S}$,
$\Ex{t(X)^2}\leq \norm{\bayes}_\infty\norm{t}^2\leq \norm{\bayes}_\infty\enspace.$

\subsection{Proof of \propref{AssApproxKern}}
First, since $\norm{K}_1\ge 1$
\begin{align*}
 \norm{\bayes_{k_{K,h}}}^2&=\int_{\Xbb}\paren{\int_{\Xbb}\bayes(y)\frac1{h}K\paren{\frac{x-y}{h}}dy}^2dx\\
& =\int_{\Xbb}\paren{\int_{\Xbb}\bayes(x+hz)K\paren{z}dz}^2dx\\
&\le\norm{K}_1^2 \int_{\Xbb}\paren{\int_{\Xbb}\bayes(x+hz)\frac{\absj{K\paren{z}}}{\norm{K}_1}dz}^2dx\\
&\le \norm{K}_1^2\int_{\Xbb^2}\bayes(x+hz)^2\frac{\absj{K\paren{z}}}{\norm{K}_1}dxdz\le \norm{\bayes}_\infty\norm{K}_1^2\enspace.
\end{align*}
Hence, Assumption~\eqref{eq:AssOr1} holds if $\Upsilon\ge 1+\norm{s}_\infty\norm{K}_1^2$. Now, we have
\[P\paren{\chi_{k_{K,h}}^2}=\frac{K(0)^2}{h^2}=P\Theta_{k_{K,h}}\frac{K(0)^2}{\norm{K}^2 h}\le nP\Theta_{k_{K,h}}\frac{K(0)^2}{\norm{K}^2 \norm{K}_\infty\norm{K}_1} \enspace,\]
so it is sufficient to have $\Upsilon\ge K(0)/\norm{K}^2$  (since $K(0)\leq \norm{K}_\infty$) to ensure \eqref{eq:AssOr2}. Moreover, for any $h\in \sH$ and any $x\in \Xbb$,
\[\bayes_{k_{K,h}}(x)=\int_{\Xbb}\bayes(y)\frac1{h}K\paren{\frac{x-y}{h}}dy=\int_{\Xbb}\bayes(x+zh)K(z)dz\leq \norm{\bayes}_\infty\norm{K}_1 \enspace.
\]
Therefore, Assumption~\eqref{eq:AssOr3} holds for $\Upsilon\ge 2\norm{\bayes}_\infty\norm{K}_1$. Then on one hand
\begin{align*}
\absj{A_{k_{K,h}}(x,y)}&\le \frac1{h^2}\int_{\Xbb}\absj{K\paren{\frac{x-z}{h}}K\paren{\frac{y-z}{h}}}dz\\
&\le \frac1{h}\int_{\Xbb}\absj{K\paren{\frac{x-y}{h}-u}K\paren{u}}du \\
&\le\frac{\norm{K}^2}h\wedge \frac{\norm{K}_\infty\norm{K}_1}{h}\le P\Theta_{k_{K,h}}\wedge n\enspace.
\end{align*}
And on the other hand
\begin{align*}
\Ex{\absj{A_{k_{K,h}}(X,x)}}&\le \frac1{h}\int_{\Xbb^2}\absj{K\paren{\frac{x-y}{h}-u}K\paren{u}}du~\bayes(y)dy\\
&=\int_{\Xbb^2}\absj{K\paren{v}K\paren{u}}\bayes(x+h(v-u))du dv\le \norm{\bayes}_\infty\norm{K}_1^2 \enspace.
\end{align*}
Therefore,
\begin{multline*}
\sup_{x\in\Xbb}~\Ex{A_{k_{K,h}}(X,x)^2}\le \sup_{(x,y)\in\Xbb^2}\absj{A_{k_{K,h}}(x,y)}~\sup_{x\in\Xbb}~\Ex{\absj{A_{k_{K,h}}(X,x)}}\\
\le \paren{P\Theta_{k_{K,h}}\wedge n}\norm{\bayes}_\infty\norm{K}_1^2\enspace,
\end{multline*}
and
$
\Ex{A_{k_{K,h}}(X,Y)^2}\le \sup_{x\in\Xbb}~\Ex{A_{k_{K,h}}(X,x)^2}\le \norm{\bayes}_\infty\norm{K}_1^2P\Theta_{k_{K,h}} \enspace. 
$
Hence Assumption~\eqref{eq:AssOr4} and \eqref{eq:AssOr5} hold when $\Upsilon \ge  \norm{\bayes}_\infty\norm{K}_1^2$. Finally let us prove that Assumption~\eqref{eq:AssVar} is satisfied. Let $t\in\Boule_{k_{K,h}}$ and $a\in L^2$ be such that $\norm{a}=1$ and $t(y)=\int_{\Xbb}a(x)\frac1{h}K\paren{\frac{x-y}{h}}dx$ for all $y\in \Xbb$. Then the following follows from Cauchy-Schwarz inequality 
\[ t(y)\leq \frac1{h}\sqrt{\int_{\Xbb}a(x)^2dx}\sqrt{\int_{\Xbb}K\paren{\frac{x-y}{h}}^2dx}\leq \frac{\norm{K}}{\sqrt{h}}\enspace. 
\]
Thus for any $t\in \Boule_{k_{K,h}}$
\[Pt^2\leq \norm{t}_\infty\psh{\absj{t}}{\bayes}\le \frac{\norm{K}}{\sqrt{h}} \norm{\bayes}= \norm{\bayes}\sqrt{P\Theta_{k_{K,h}}}\le \sqrt{\Upsilon P\Theta_{k_{K,h}}}\enspace.\]
We conclude that all the assumptions hold if 
\[\Upsilon \ge \paren{K(0)/\norm{K}^2}\vee \paren{1+2\norm{\bayes}_\infty\norm{K}_1^2}\enspace.\]

\subsection{Proof of \propref{AssWeiProjKern}}
Let us define for convenience $\Phi(x)\egaldef\sum_{i=1}^p\vphi_i(x)^2$, so $\Gamma \ge 1\vee n^{-1} \norm{\Phi}_\infty$. Then we have for these kernels: $ \Phi(x)\ge \chi_{k_w}(x) \ge\Theta_{k_w}(x)$ for all $x\in\Xbb$. Moreover, denoting by $\Pi\bayes$ the orthogonal projection of $\bayes$ onto the linear span of $(\vphi_i)_{i=1,\ldots,p}$,
\begin{align*}
 \norm{\bayes_{k_w}}^2=\sum_{i=1}^pw_i^2\paren{P\vphi_i}^2\le \norm{\Pi\bayes}^2\le \norm{\bayes}^2\le \norm{\bayes}_\infty\enspace.
\end{align*}
Assumption~\eqref{eq:AssOr1} holds for this family if $\Upsilon\ge \Gamma(1+\norm{s}_\infty)$. We prove in what follows that all the remaining assumptions are valid using only \eqref{eq:CdtForGamma} and \eqref{eq:AssOr1}.\\
First, it follows from Cauchy-Schwarz inequality that, for any $x\in\Xbb$, $\chi_{k_w}(x)^2\le  \Phi(x) \Theta_{k_w}(x)$. Assumption~\eqref{eq:AssOr2} is then automatically satisfied from the definition of $\Gamma$ 
\[ \Ex{\chi_{k_w}(X)^2}\le \norm{ \Phi}_\infty P\Theta_{k_w} \le \Gamma n P\Theta_{k_w} \enspace. \] 
\noindent
Now let $w$ and $w'$ be any two vectors in $[0,1]^p$, we have
\[\bayes_{k_w}=\sum_{i=1}^pw_i(P\vphi_i)\vphi_i,\qquad \bayes_{k_w}-\bayes_{k_{w'}}=\sum_{i=1}^p(w_i-w_i')\paren{P\vphi_i}\vphi_i\enspace.\]
Hence
$\norm{\bayes_{k_w}-\bayes_{k_{w'}}}^2=\sum_{i=1}^p(w_i-w_i')^2\paren{P\vphi_i}^2$
and, by Cauchy-Schwarz inequality, for any $x\in \Xbb$,
\[\absj{\bayes_{k_w}(x)-\bayes_{k_{w'}}(x)}\le \norm{\bayes_{k_w}-\bayes_{k_{w'}}}\sqrt{\Phi(x)}\le \norm{\bayes_{k_w}-\bayes_{k_{w'}}}\sqrt{\Gamma n}\enspace.\]
Assumption~\eqref{eq:AssOr3} follows using \eqref{eq:AssOr1}. Concerning Assumptions~\eqref{eq:AssOr4} and \eqref{eq:AssOr5}, let us first notice that by orthonormality, for any $(x,x')\in\Xbb^2$,
\vspace{-0.1cm}
\[A_{k_w}(x,x')=\sum_{i=1}^pw_i^2\vphi_i(x)\vphi_i(x')\enspace. \]
\vspace{-0.3cm}

\noindent Therefore, Assumption~\eqref{eq:AssOr5} holds since
\vspace{-0.2cm}
\begin{align*}
\Ex{A_{k_w}(X,x)^2}&=\int_\Xbb\paren{\sum_{i=1}^pw_i^2\vphi_i(y)\vphi_i(x)}^2\bayes(y)d\mu(y)\\
&\le \norm{\bayes}_\infty \sum_{1\le i,j\le p}w_i^2w_j^2\vphi_i(x)\vphi_j(x)\int_{\Xbb}\vphi_i(y)\vphi_j(y)d\mu(y)\\
&= \norm{\bayes}_\infty \sum_{i=1}^p w_i^4\vphi_i(x)^2\le \norm{\bayes}_\infty \Phi(x)\le \norm{\bayes}_\infty \Gamma n \enspace.
\end{align*}
\vspace{-0.4cm}

\noindent Assumption~\eqref{eq:AssOr4} also holds from similar computations:
\vspace{-0.25cm}
\begin{align*}
 \Ex{A_{k_w}(X,Y)^2}&=\int_\Xbb\Ex{\paren{\sum_{i=1}^pw_i^2\vphi_i(X)\vphi_i(x)}^2}\bayes(x)d\mu(x)\\
 &\le \norm{\bayes}_\infty\sum_{1\le i,j\le p}w_i^2w_j^2\Ex{\vphi_i(X)\vphi_j(X)}\int_{\Xbb}\vphi_i(x)\vphi_j(x)d\mu(x)\\
 &\le \norm{\bayes}_\infty P\Theta_{k_w}\enspace.
\end{align*}
\vspace{-0.4cm}

\noindent
We finish with the proof of \eqref{eq:AssVar}. Let us prove that $\Boule_{k_w}=\sE_{k_w}$, where
\vspace{-0.15cm}
\begin{equation*}
\sE_{k_w}=\set{t=\sum_{i=1}^pw_it_i\vphi_i,\telque \sum_{i=1}^pt_i^2\le 1}\enspace. 
\end{equation*}
\vspace{-0.4cm}

\noindent
First, notice that any $t\in  \Boule_{k_w}$ can be written
\vspace{-0.15cm}
\begin{align*}
 \int_{\Xbb}a(x)k_w(x,y)d\mu(x)=\sum_{i=1}^pw_i\psh{a}{\vphi_i}\vphi_i(y)\enspace.
\end{align*}
\vspace{-0.4cm}

\noindent
Then, consider some $t\in \sE_{k_w}$. By definition, there exists a collection $(t_i)_{i=1,\ldots,p}$ such that $t=\sum_{i=1}^pw_it_i\vphi_i$, and $\sum_{i=1}^pt_i^2\le 1$. If $a=\sum_{i=1}^pt_i\vphi_i$, $\norm{a}^2= \sum_{i=1}^pt_i^2 \le 1$ and $\psh{a}{\vphi_i}=t_i$, hence $t\in \Boule_{k_w}$.
Conversely, for $t\in \Boule_{k_w}$, there exists some function $a\in L^2$ such that $\norm{a}^2\le 1$, and $t=\sum_{i=1}^pw_i\psh{a}{\vphi_i}\vphi_i$. Since $(\vphi_i)_{i=1,\ldots,p}$ is an orthonormal system, one can take $a=\sum_{i=1}^p \psh{a}{\vphi_i} \vphi_i$. With $t_i=\psh{a}{\vphi_i}$, we find $\norm{a}^2= \sum_{i=1}^pt_i^2$ and $t\in \sE_{k_w}$.
For any $t\in\Boule_{k_w}=\sE_{k_w}$, $\norm{t}^2=\sum_{i=1}^pw_i^2 t_i^2\le \sum_{i=1}^p t_i^2\le 1$. Hence
$Pt^2\leq \norm{\bayes}_\infty\norm{t}^2\le \norm{\bayes}_{\infty}\enspace.$

\vspace{-0.5cm}
\section{Concentration of the residual terms}
The following proposition gathers the concentration bounds of the remaining terms appearing in \eqref{eq:DecRisk}. 
\begin{proposition}\label{prop:ConcRemainder}
Let $\set{k}_{\kK}$ denote a finite collection of kernels satisfying \eqref{eq:CdtForGamma} and suppose that Assumptions~\eqref{eq:AssOr1}, \eqref{eq:AssOr2} and \eqref{eq:AssOr3} hold. Then
 \begin{equation}\label{eq:Psk}
 \forall \theta\in(0,1),\qquad 2\frac{P(\bayes_{\kh}-\bayes_k)}n\le \theta \norm{\bayes-\bayes_{\kh}}^2+\theta \norm{\bayes-\bayes_{k}}^2+\frac{2\Upsilon}{\theta n}\enspace.
 \end{equation}
  For any $x\ge 1$, with probability larger than $1-2\absj{\sK}^2e^{-x}$, for any $(k,k')\in\sK^2$, for any $\theta\in(0,1)$,
 \begin{equation}\label{eq:Pn-Psk}
 \absj{2(P_n-P)(\bayes_k-\bayes_{k'})}\le \theta\paren{ \norm{\bayes-\bayes_{k'}}^2+\norm{\bayes-\bayes_{k}}^2}+\square \frac{\Upsilon x^2}{ \theta n}\enspace.
 \end{equation}
 For any $x\ge 1$, with probability larger than $1-2\absj{\sK}e^{-x}$, for any $\kK$,
  \begin{equation}\label{eq:Pn-Pchik}
\forall\theta\in(0,1),\qquad \absj{2(P_n-P)\chi_k}\le \theta P\Theta_k+\square\frac{\Upsilon x}{\theta}\enspace.
 \end{equation}
 For any $x\ge 1$, with probability larger than $1-5.4\absj{\sK}e^{-x}$, for any $\kK$, 
   \begin{equation}\label{eq:Uk}
\forall\theta\in(0,1),\qquad\frac{2\absj{U_k}}{n^2}\le \theta\frac{P\Theta_k}n+ \square \frac{\Upsilon x^2}{\theta n} \enspace.
\end{equation}
\end{proposition}
{\it Proof}
First for \eqref{eq:Psk}, notice that, by \eqref{eq:AssOr3}, for any $\theta\in(0,1)$
\begin{align*}
2\frac{P(\bayes_{\kh}-\bayes_k)}n&\le2\frac{\norm{\bayes_{\kh}-\bayes_k}_\infty}n \le\frac2n\paren{\Upsilon\vee\paren{\frac{\theta}4 n\norm{\bayes_k-\bayes_{\kh}}^2+\frac{\Upsilon}{\theta}}} \\
&\le  \frac{\theta}2 \norm{\bayes_k-\bayes_{\kh}}^2+\frac{2\Upsilon}{\theta n}\le \theta \norm{\bayes-\bayes_{\kh}}^2+\theta \norm{\bayes-\bayes_{k}}^2+\frac{2\Upsilon}{\theta n}\enspace.
\end{align*}
Then, by \propref{Bern}, with probability larger than $1-\absj{\sK}^2e^{-x}$, 
\[\mbox{for any }(k,k')\in\sK^2,\quad (P_n-P)(\bayes_k-\bayes_{k'})\le \sqrt{\frac{2P\paren{\bayes_k-\bayes_{k'}}^2 x}n}+\frac{\norm{\bayes_k-\bayes_{k'}}_\infty x}{3n}\enspace.\]
Since by \eqref{eq:AssOr1}
$ P\paren{\bayes_k-\bayes_{k'}}^2\le \norm{\bayes}_\infty\norm{\bayes_k-\bayes_{k'}}^2\le \Upsilon\norm{\bayes_k-\bayes_{k'}}^2\enspace,$
\[\sqrt{\frac{2P\paren{\bayes_k-\bayes_{k'}}^2x}n}\le \frac{\theta}4\norm{\bayes_k-\bayes_{k'}}^2+\frac{2\Upsilon x}{\theta n}\enspace.\]
Moreover, by \eqref{eq:AssOr3}
$\frac{\norm{\bayes_k-\bayes_{k'}}_\infty x}{3n}\le \frac{\theta}4 \norm{\bayes_k-\bayes_{k'}}^2+\square\frac{\Upsilon x^2}{ \theta n}\enspace.$
Hence, for $x\ge 1$, with probability larger than $1-\absj{\sK}^2e^{-x}$
\vspace{-0.2cm}
\begin{align*}
(P_n-P)(\bayes_k-\bayes_{k'})&\le \frac{\theta}2 \norm{\bayes_k-\bayes_{k'}}^2+\square\frac{\Upsilon x^2}{ \theta n}\\
&\le \theta\paren{ \norm{\bayes-\bayes_{k'}}^2+\norm{\bayes-\bayes_{k}}^2}+\square\frac{\Upsilon x^2}{ \theta n}\enspace, 
\end{align*}
\vspace{-0.4cm}

\noindent which gives \eqref{eq:Pn-Psk}. 
Now, using again \propref{Bern}, with probability larger than $1-\absj{\sK}e^{-x}$, for any $\kK$,
\[(P_n-P)\chi_k\le \sqrt{\frac{2P\paren{\chi_k}^2 x}n}+\frac{\norm{\chi_k}_\infty x}{3n}\enspace.\]
By \eqref{eq:CdtForGamma} and \eqref{eq:AssOr1},  for any $\kK$,
$\norm{\chi_k}_\infty\le \sup_{(x,y)\in\Xbb^2}\absj{k(x,y)}\le \Gamma n\le \Upsilon n\enspace.$

\noindent Concerning \eqref{eq:Pn-Pchik}, we get by \eqref{eq:AssOr2}, $P\chi_k^2\le \Upsilon n P\Theta_k$,
hence, for any $x\ge 1$ we have with probability larger than $1-\absj{\sK}e^{-x}$
\[(P_n-P)\chi_k\le \theta P\Theta_k+\paren{\frac13+\frac{1}{2\theta}}\Upsilon x\enspace.\]
For \eqref{eq:Uk}, we apply \propref{ConcUstat} to obtain with probability larger than $1-2.7\absj{\sK}e^{-x}$, for any $\kK$,
\begin{align*}
 \frac{U_k}{n^2}\le \frac{\square}{n^2}\paren{C\sqrt{x}+Dx+Bx^{3/2}+Ax^2}\enspace,
\end{align*}
where $A,B,C,D$ are defined accordingly to \propref{ConcUstat}.  Let us evaluate all these terms.
First, $A\le 4\sup_{(x,y)\in\Xbb^2}\absj{k(x,y)}\le 4\Upsilon n$ by \eqref{eq:CdtForGamma} and \eqref{eq:AssOr1}. Next,
$
 C^2\le \square n^2\Ex{k(X,Y)^2}\le \square n^2\norm{\bayes}_\infty P\Theta_k\le \square n^2\Upsilon P\Theta_k\enspace.
$

\noindent Using \eqref{eq:CdtForGamma}, we find 
$B^2\leq 4 n \sup_{x\in \Xbb} \int k(x,y)^2 \bayes(y) d\mu(y) \leq 4 n\norm{\bayes}_\infty \Gamma \enspace.$

\noindent By \eqref{eq:AssOr1}, we consequently have $ B^2\le 4 \Upsilon n$.
Finally, using Cauchy-Schwarz inequality and proceeding as for $C^2$,
\[\Ex{\sum_{i=1}^{n-1}\sum_{j=i+1}^na_i(X_i)b_j(X_j)k(X_i,X_j)}\le n\sqrt{\Ex{k(X,Y)^2}}\le n\sqrt{\Upsilon P\Theta_k}\enspace.\]
Hence, $D\le n\sqrt{\Upsilon P\Theta_k}$ which gives \eqref{eq:Uk}.


\begin{thebibliography}

\bibitem[AB09]{Arl_Bac:2009}
S.~Arlot and F.~Bach.
\newblock Data-driven calibration of linear estimators with minimal penalties.
\newblock In {\em Advances in Neural Information Processing Systems 22}, pages
  46--54, 2009.

\bibitem[Ada06]{Adamczak2006}
R.~Adamczak.
\newblock Moment inequalities for {$U$}-statistics.
\newblock {\em Ann. Probab.}, 34(6):2288--2314, 2006.

\bibitem[AM09]{Arl_Mas:2009}
S.~Arlot and P.~Massart.
\newblock Data-driven calibration of penalties for least-squares regression.
\newblock {\em J. Mach. Learn. Res.}, 10:245--279, 2009.

\bibitem[Bir06]{Birge2006}
L.~Birg{\'e}.
\newblock Statistical estimation with model selection.
\newblock {\em Indag. Math. (N.S.)}, 17(4):497--537, 2006.

\bibitem[BLM13]{BLM2013}
S.~Boucheron, G.~Lugosi, and P.~Massart.
\newblock {\em Concentration inequalities}.
\newblock Oxford University Press, Oxford, 2013.

\bibitem[BM07]{Bir_Mas:2007}
L.~Birg{\'e} and P.~Massart.
\newblock Minimal penalties for {G}aussian model selection.
\newblock {\em Probab. Theory Related Fields}, 138(1-2):33--73, 2007.

\bibitem[DJKP96]{DJKP1996}
D.~L. Donoho, I.~M. Johnstone, G.~Kerkyacharian, and D.~Picard.
\newblock Density estimation by wavelet thresholding.
\newblock {\em Ann. Statist.}, 24(2):508--539, 1996.

\bibitem[DL01]{Devroye-Lugosi2001}
L.~Devroye and G.~Lugosi.
\newblock {\em Combinatorial methods in density estimation}.
\newblock Springer Series in Statistics. Springer-Verlag, 2001.

\bibitem[DO13]{DehOuadah2013}
P.~Deheuvels and S.~Ouadah.
\newblock Uniform-in-bandwidth functional limit laws.
\newblock {\em J. Theoret. Probab.}, 26(3):697--721, 2013.

\bibitem[EL99a]{Egg_LaR:1999b}
P.~P.~B. Eggermont and V.~N. LaRiccia.
\newblock Best asymptotic normality of the kernel density entropy estimator for
  smooth densities.
\newblock {\em IEEE Trans. Inform. Theory}, 45(4):1321--1326, 1999.

\bibitem[EL99b]{Egg_LaR:1999}
P.~P.~B. Eggermont and V.~N. LaRiccia.
\newblock Optimal convergence rates for good's nonparametric maximum likelihood
  density estimator.
\newblock {\em Ann. Statist.}, 27(5):1600--1615, 1999.

\bibitem[EL01]{Egg_LaR:2001}
P.~P.~B. Eggermont and V.~N. LaRiccia.
\newblock {\em Maximum penalized likelihood estimation}, volume~I of {\em
  Springer Series in Statistics}.
\newblock Springer-Verlag, New York, 2001.

\bibitem[FT06]{FrTul2006}
M.~Fromont and C.~Tuleau.
\newblock Functional classification with margin conditions.
\newblock In {\em Learning theory}, volume 4005 of {\em Lecture Notes in
  Comput. Sci.}, pages 94--108. Springer, Berlin, 2006.

\bibitem[GL11]{GL2011}
A.~Goldenshluger and O.~Lepski.
\newblock Bandwidth selection in kernel density estimation: oracle inequalities
  and adaptive minimax optimality.
\newblock {\em Ann. Statist.}, 39(3):1608--1632, 2011.

\bibitem[GLZ00]{GLZ2000}
E.~Gin{\'e}, R.~Lata{\l}a, and J.~Zinn.
\newblock Exponential and moment inequalities for {$U$}-statistics.
\newblock In {\em High dimensional probability, {II}}, volume~47 of {\em Progr.
  Probab.}, pages 13--38. Birkh\"auser Boston, 2000.

\bibitem[GN09]{GinNickl2009}
E.~Gin{\'e} and R.~Nickl.
\newblock Uniform limit theorems for wavelet density estimators.
\newblock {\em Ann. Probab.}, 37(4):1605--1646, 2009.

\bibitem[GN15]{Gin_Nic:2015}
E~Gin{\'e} and R.~Nickl.
\newblock {\em Mathematical foundations of infinite-dimensional statistical
  models}.
\newblock Cambridge University press, 2015.

\bibitem[HRB03]{HRB2003}
C.~Houdr{\'e} and P.~Reynaud-Bouret.
\newblock Exponential inequalities, with constants, for {U}-statistics of order
  two.
\newblock In {\em Stochastic inequalities and applications}, volume~56 of {\em
  Progr. Probab.}, pages 55--69. Birkh\"auser, Basel, 2003.

\bibitem[Ler12]{Le2012}
M.~Lerasle.
\newblock Optimal model selection in density estimation.
\newblock {\em Ann. Inst. H. Poincar{\'e} Probab. Statist.}, 48(3):884--908,
  2012.

\bibitem[Mas07]{Massart2007}
P.~Massart.
\newblock {\em Concentration inequalities and model selection}, volume 1896 of
  {\em Lecture Notes in Mathematics}.
\newblock Springer, Berlin, 2007.
\newblock Lectures from the 33rd Summer School in Saint-Flour.

\bibitem[MS11]{MasSwan2011}
D.~M. Mason and J.~W.~H. Swanepoel.
\newblock A general result on the uniform in bandwidth consistency of
  kernel-type function estimators.
\newblock {\em TEST}, 20(1):72--94, 2011.

\bibitem[MS15]{MasSwan2015}
D.~M. Mason and J.~W.~H. Swanepoel.
\newblock Erratum to: A general result on the uniform in bandwidth consistency
  of kernel-type function estimators .
\newblock {\em TEST}, 24(1):205--206, 2015.

\bibitem[Rig06]{Rig2006}
P.~Rigollet.
\newblock Adaptive density estimation using the blockwise {S}tein method.
\newblock {\em Bernoulli}, 12(2):351--370, 2006.

\bibitem[RT07]{RigTsy2007}
P.~Rigollet and A.~B. Tsybakov.
\newblock Linear and convex aggregation of density estimators.
\newblock {\em Math. Methods Statist.}, 16(3):260--280, 2007.

\bibitem[Tsy09]{Tsy:2009}
A.~B. Tsybakov.
\newblock {\em Introduction to nonparametric estimation}.
\newblock Springer Series in Statistics. Springer, New York, 2009.

\end{thebibliography}
\end{document}